\documentclass{amsart}

\usepackage{amsmath, amssymb}
\usepackage{graphicx}
\usepackage{tikz} 
\usepackage{multirow}
\usepackage{array,longtable}
\usetikzlibrary{matrix,arrows}
\usepackage{rotating}
\usepackage{float}

\newtheorem{theorem}{Theorem}[section]
\newtheorem{lemma}[theorem]{Lemma}
\newtheorem{corollary}[theorem]{Corollary}
\newtheorem{proposition}[theorem]{Proposition}
\newtheorem{conjecture}[theorem]{Conjecture}
\theoremstyle{definition}
\newtheorem{definition}[theorem]{Definition}
\newtheorem{example}[theorem]{Example}

\theoremstyle{remark}
\newtheorem{remark}[theorem]{Remark}
\newtheorem{notation}[theorem]{Notation}
\numberwithin{equation}{section}

\newcommand{\Or}{\mathcal O}
\newcommand{\R}{\mathbb R}
\newcommand{\Z}{\mathbb Z}

\newcommand{\gc}{\mathfrak{g}}

\newcommand{\mf}{\mathfrak}

\newcommand{\tu}{\widetilde}

\newcommand{\Br}{\overline}

\newcommand{\la}{\langle}
\newcommand{\ra}{\rangle}
\newcommand{\SH}{ \mathcal{SH}}
\newcommand{\mca}{\mathcal}
\newcommand{\Smal}{\prod_{\rho/2} ^s (\tu{G})} 

\newcommand{\ths}{\thinspace}
\newcommand{\ul}{\underline}
\begin{document}

\title{Some Genuine Small Representations of a Nonlinear Double Cover}

\author{Wan-Yu Tsai}
\address{Institute of Mathematics, Academia Sinica, Taipei, TAIWAN }
\curraddr{Institute of Mathematics, Academia Sinica, 6F, Astronomy-Mathematics Building, No. 1, Sec. 4, Roosevelt Road, Taipei 10617, TAIWAN }
\email{wytsai@math.sinica.edu.tw}

%\author{Author Two}
%\address{Mathematical Research Section, School of Mathematical Sciences,
%Australian National University, Canberra ACT 2601, Australia}

%    General info
%\subjclass[2000]{Primary 54C40, 14E20; Secondary 46E25, 20C20}

%\date{January 1, 2001 and, in revised form, June 22, 2001.}

%\dedicatory{This paper is dedicated to Jeffrey Adams}

%\keywords{Lie groups, nonlinear cover, representation theory}

\begin{abstract}
Let $G$ be the real points of a simply connected, semisimple, simply laced complex Lie group, and let $\widetilde{G}$ be the nonlinear double cover of $G$. We discuss a set of small genuine irreducible representations of $\widetilde{G}$ which  can be characterized by the following properties: (a) the infinitesimal character is $\rho/2$; (b) they have maximal $\tau$-invariant; (c) they have a particular associated variety $\mathcal{O}$. When $G$ is split, we construct them explicitly. Furthermore, in many cases,  there is a one-to-one correspondence between  these small representations and the pairs (genuine central characters of $\tu{G}$, real forms of $\mathcal{O}$) via the map $\tu{\pi} \mapsto (\chi_{\tu{\pi}}, AV(\tu{\pi}) )$.\end{abstract}

\maketitle

%\section*{This is an unnumbered first-level section head}
%This is an example of an unnumbered first-level heading.

%% The correct journal style for \specialsection is all uppercase; a known bug
%% in amsart.cls prevents this, so input must be uppercase until it is fixed.
%\specialsection*{This is a Special Section Head}
%\specialsection*{THIS IS A SPECIAL SECTION HEAD}
%This is an example of a special section head%
%%%%%%%%%%%%%%%%%%%%%%%%%%%%%%%%%%%%%%%%%%%%%%%%%%%%%%%%%%%%%%%%%%%%%%%%
%\footnote{Here is an example of a footnote. Notice that this footnote
%text is running on so that it can stand as an example of how a footnote
%with separate paragraphs should be written.
%\par
%And here is the beginning of the second paragraph.}%
%%%%%%%%%%%%%%%%%%%%%%%%%%%%%%%%%%%%%%%%%%%%%%%%%%%%%%%%%%%%%%%%%%%%%%%%
.

\section{Introduction}

Assume that $G_{\mathbb{C}}$ is a simply connected, semisimple complex Lie group, and $G$ is a real form of $G_{\mathbb{C}}$ with nontrivial fundamental group. Then $G$ has a nonlinear double cover $\widetilde{G}$, which is not a matrix group (see \cite{AT}, Proposition 3.6).  For example, this holds when $G$ is a split group, or $G=SU(p,q)$, $Spin(p,q)$, $p, q\ge 2$, and most  real forms of the exceptional  groups. In fact, most real forms of $G_{\mathbb{C}}$ have a nonlinear double cover (see \cite{A2}). The purpose of  this paper is to discuss some small genuine representations of $\widetilde{G}$ and their properties. By a genuine representation of  $\tu{G}$ we mean that a representation of $\widetilde{G}$ which does not factor through $G$.

The paper studies the set of irreducible admissible genuine representations of $\tu G$ with a fixed infinitesimal character $\lambda$ (in particular, $\lambda =\rho/2$ when $G$ is simply laced)
and having maximal $\tau$-invariants, which is the notion to define "\textit{small}" in this paper. The number of such set of representations is usually small, which can be calculated by  the coherent continuation representation carried by the integral Weyl group. The study of these small representations is important and interesting due to a number of reasons. We list some of them
as follows.

\begin{enumerate}
\item[1.] Minimal representations of reductive groups have been studied for long time and are of significant importance in many areas, such as number theory and physics (e.g. see \cite{Ka}, \cite{KPW}). They are the 
smallest infinite dimensional irreducible unitary representations. More precisely, a minimal representation $\pi$ of a reductive group $G$ is characterized by either of the following properties. 
\begin{itemize}
\item[(a)]
The annihilator $\mathcal J$ of $\pi$ in the universal enveloping algebra $U(\gc)$ is equal to the Joseph ideal (see \cite{J2}), with $\gc$ the complexified Lie algebra of $G$.
\item[(b)]  The associated variety of $\mathcal J$ is the closure of the minimal $G$-coadjoint orbit in $\gc ^*$.
\item[(c)] The Gelfand-Kirillov dimension of $\pi$ attains its minimum among all infinite dimensional representations of $G$. 
\end{itemize}
There has been an active study on minimal representations, either by algebraic approaches (e.g. \cite{GS}, \cite{J2}, \cite{Ka}, \cite{KO1}),  or by analytic approaches (e.g. \cite{HKMO}, \cite{KM}, \cite{KO2}, \cite{KO3}). 

When the group is the metaplectic double cover of $Sp(2n, \R)$, the minimal representations are the 
 Weil representations, which play a central role in the theory of automorphic forms and theta  correspondence. In this case, the minimal representations are precisely the small representations described in this paper. 
 For other groups, the small representations we study are not attached to  the minimal nilpotent orbit but to one of the next smallest possibilities. As with the minimal representations, they can occur as local
 constituents of metaplectic automorphic representations, and there will be applications. 
 There are some works with applications of non-minimal small representations, arising from a nonlinear cover of some other linear groups 
(see e.g. \cite{BFG}, \cite{KP1}, \cite{LS1}, \cite{T1}).

\item[2.] For split groups, part of the genuine small representations are the unique irreducible quotients of pseudospherical principal series representations which are the Shimura lifts of the trivial representation of the 
dual split group $G^{\vee} (\R)$ (see \cite{ABPTV}). 
 
\item[3.] It is expected that most of these small representations have multiplicity free $K$-types. See \cite{Luc} and Tables 9, 10, 11 for the explicit $K$-structure formulas for the split groups of type $A$ and $D$. Indeed, the  representations  annihilated by the maximal primitive  ideal with infinitesimal character $\rho/2$ are expected to have $K$-multiplicity free.  For example, a precise conjecture (Conjecture  8.14 in \cite{AHV}) is stated for the split group of type $E_8$. 
For split groups in general, see \cite{LS2}. The $K$-structure of such a small representation $\pi$ should be almost equal to the $K$-representation carried by the regular functions on the associated variety of $\pi$. For type $D_n$, this is demonstrated in recent preprints \cite{BTs1} and \cite{BTs2}.

\item[4.] Just as minimal representations are one of the building blocks of unitary representations of Lie groups,  small representations are expected to be 
unitary, and they should play a role in the understanding of the unitary dual of $\tu G$. For instance,  see \cite{ABPTV}, \cite{Hu}, \cite{Luc}, and \cite{LS1} for some cases of classical groups.
We hope to come back to prove the unitarity of the small representations described in this paper in the future. 
 
\end{enumerate}

Here we give an outline of the paper.

    In  Section 2, we first introduce some basic invariants, such as infinitesimal character, $\tau$-invariant, and associated variety, which are  used to classify representations. These notions are quite general, and are defined for more general real reductive groups $G$, which can be linear or nonlinear. For each type, we fix an infinitesimal character $\lambda$ as listed in Table 1. 
These infinitesimal characters are chosen such that the irreducible quotients of the pseudospherical principal series representations defined in \cite{ABPTV} are \textit{small}.     
Then we define a class of representations of $\widetilde{G}$ denoted
  
    \begin{center} $\prod _{\lambda} ^s (\widetilde{G})=\{\widetilde{  \pi} \thinspace |   \thinspace \widetilde{\pi }\in \widehat{\widetilde{G}}_{adm,\lambda}, \text{ $\widetilde{\pi}$ is genuine and has maximal $\tau$-invariant} \}$, \end{center} where $\widehat{\widetilde{G}}_{adm,\lambda}$ is the set of irreducible admissible representations of $\widetilde{G}$ with infinitesimal character $\lambda$. 
Here the superscript $s$ stands for \textit{small} in the sense that the representations in this set have maximal $\tau$-invariant. There is a unique complex nilpotent orbit $\mathcal{O}$ which is the complex associated variety of  every $\widetilde{\pi}$ from $\prod _{\lambda} ^s (\widetilde{G})$, i.e. $AV(I_{\tu{\pi}})=\Br{ \mathcal{O} }$ (see the notation in Section \ref{AVGKdim}). We calculate this orbit $\mathcal{O}$  explicitly  for all types and  list them in Table 1. 
    
        Denote\begin{center} $\prod _{\lambda} ^{\mathcal{O}} (\widetilde{G})= \{\widetilde{\pi} \thinspace |   \thinspace\widetilde{ \pi }\in \widehat{\widetilde{G}}_{adm,\lambda}, \widetilde{\pi}\text{ is genuine and  $AV(I_{\tu{\pi}})=\Br{ \mathcal{O} }$} \}$. \end{center} 

      Then we have
\begin{theorem} \label{maxtauO} 
$\prod_{\lambda} ^s (\widetilde{G})=\prod_{\lambda} ^{\mathcal{O}} (\widetilde{G})$.
\end{theorem}

The proof of the theorem is based  on truncated induction of representations of Weyl groups and the Springer correspondence. 

This set of representations $\prod_{\lambda} ^s (\widetilde{G})=\prod_{\lambda} ^{\mathcal{O}} (\widetilde{G})$ is what we are going to discuss throughout this paper. First of all, we can attach to each $\widetilde{\pi} \in \prod_{\lambda} ^s (\widetilde{G})$ a pair $(\chi_{\widetilde{\pi}}, \mathcal{O}_{\widetilde{\pi}}   )$, where $\chi_{\widetilde{\pi}}$ is the central character of $\widetilde{\pi}$ and $ \mathcal{O}_{\widetilde{\pi}}$ is the real associated variety of $\widetilde{\pi}$, denoted AV$(\tu{\pi}) = \Br{\Or _{\tu{\pi}}}$. Here, $ \mathcal{O}_{\widetilde{\pi}}$ is one of the real forms of $\mathcal{O}$, and in Section 4, we will see that there are not many real groups which have nonempty intersection with $\mathcal{O}$ and the number of real forms of $\mathcal{O}$ is tiny as well.  The notions of real associated variety and genuine central character will be discussed in more detail in Sections 4 and 5. 

   In Section 6, we restrict our attention to simply laced split groups, and hence $\lambda =\rho/2$.  For split groups, there is a well-understood family of representations described in \cite{ABPTV}, namely, pseudospherical principal series representations. We consider the unique irreducible quotients of them with infinitesimal character $\rho/2$, and denote the set of these representations by $\SH _{\rho/2}$. 
    Starting with $\SH _{\rho/2}$, we construct other genuine representations in $\prod_{\rho/2} ^s (\widetilde{G})$. There are standard ways to get new representations from old ones: the theory of cross actions and Cayley transforms. In our setting these are non-standard, because they involve half-integral roots. It is possible to start with a representation from $\SH_{\rho/2}$, and apply some cross actions and Cayley transforms to it, to obtain other representations in $\prod _{\rho/2} ^s(\widetilde{G})$. The conditions which need to be satisfied are very rigid, and we get a small number of representations in $\prod _{\rho/2} ^s(\widetilde{G})$. Let $\prod _{R_D} (\widetilde{G})$ denote the set of representations obtained this way. In the last part of  Section 6, furthermore, by counting the elements in $\prod _{\rho/2} ^s(\widetilde{G})$ using  a Weyl group calculation, we show that $\Smal = \prod _{R_D} (\widetilde{G})$ for type $A_{n-1}$ and $D_n$. We conjecture this is true for type $E$.
   
   In Section 7, we denote \begin{center} $P_{\Or} (\tu{G})=\{  (\chi _i, \mathcal{O} _j ) \thinspace| \thinspace   \chi _i \in \prod _g (Z(\tu{G})), \mathcal{O}_j \text{ is a real form of } \mathcal{O} \}$,\end{center} where $\prod_g ( Z( \tu{G} ))$ is the set of genuine central characters of $\tu{G}$. The fact is that $|\Smal|=|P_{\Or} (\tu{G })|$ if $G$ has type $A_{n-1}$ or $D_n$. Furthermore, the map $\widetilde{\pi} \mapsto (\chi_{ \widetilde{\pi}  } ,\text{AV} (\tu{\pi }))$ gives a bijection between $\Smal$ and $P_{\Or}(\tu{G} )$  in many cases. Here is the main theorem.
   
\begin{theorem}   
  Let $G$ be the split real form of a simply connected,  semisimple, complex Lie group and let $\tu{G}$ be the nonlinear double cover of $G$. Consider the map \begin{center}$\xi: \widetilde{\pi} \mapsto (\chi_{ \widetilde{\pi}  } ,\text{AV} (\tu{\pi }))$\end{center} from $\Smal$ to $P_{\Or} (\tu{G})$. Then, (a) $\xi$ is bijective if $G$ has type $A_{n-1}$, unless $n$ is divisible by 4; (b) $\xi$ is bijective if $G$ has type $D_{2m}$; (c) $\xi$ is surjective if $G$ has type $E_6$ or $E_8$. 
 \end{theorem}   
   We conjecture that $\xi$ is bijective for type $E$. 
   
This paper is part of my Ph.D. thesis.  In my thesis, we introduce the Kazhdan-Patterson lifting, which is an operator taking stable representations of $G$ to 0 or virtual genuine representations of $\tu{G}$ (see \cite{AHe}).
The existence of the lifting theory was first conjectured by Kazhdan and Patterson in \cite{KP2} for the case of $GL(n, \mathbb F)$. This lifting operator is defined  to take representations of $GL(n,\mathbb F)$ to 
representations of the nonlinear covering groups of $GL(n,\mathbb F)$. In the end of Section 6 in \cite{KP2}, the conjecture says that if $\pi$ is an irreducible representation, then lifting of $\pi$ is an irreducible representation,
up to sign, or is zero. We generalize this to various real groups. 
 We claim that when $G$ is simply laced and split, the small representations in $\Smal$ can be obtained from lifting of the trivial representation (see \cite{AHu} for the result of the case when $G=GL(n,\R)$). 
 Indeed, the representations that we study in the paper for type $A$ appear in \cite{KP1} as the local component  of the generalized theta series (see Theorem I.6.4 (b) in \cite{KP1}). 
 Since the lifting operator is worked on the level of global characters, we compute the character formulas for these small representations, and expect to prove the unitarity of these representations as in the case of $GL(n,\R)$. This part still in progress is to appear in a future paper.     
 
\textit{Acknowledgement.} The author is grateful to her advisor, Professor Jeffrey Adams, for his guidance and patience, as well as Dan Barbasch and Peter Trapa for helpful discussions.

\section{Invariants of a representation} 
 
 Before introducing the set of representations of interest, some notions are needed. Let us get started with the setting. Let $G$ be a connected  real Lie group, and suppose that the complexified Lie algebra of $G$, denoted $\mathfrak{g}$, is reductive. Here  $G$ is allowed to be nonlinear, which means it cannot be embedded  into any $GL(n, \mathbb{C})$ (see \cite{AT}, for example). We fix a Cartan involution $\theta$ of $G$ and let $K=G^{\theta}$ be the corresponding maximal compact subgroup. Let $\mathfrak{h}$ be a Cartan subalgebra of $\mathfrak{g}$, and $U(\mathfrak{g})$ be the enveloping algebra of $\mathfrak{g}$. Let $\triangle=\triangle(\mathfrak{g},\mathfrak{h})$ be the root system and $W$ be the Weyl group of $\mathfrak{g}$. 

  Let $\mathcal{HC}(\mathfrak{g}, K)$ be the set of Harish-Chandra modules and let $\widehat{G}_{adm}$ denote the set of equivalence classes of irreducible admissible representations of $G$. Then   $\widehat{G}_{adm}$ can be viewed as a subset of $\mathcal{HC}(\mathfrak{g}, K)$  by sending an irreducible admissible representation $\pi\in \widehat{G}_{adm}$ to its space $V_{\pi}$ of $K$-finite vectors and then the latter can be regarded as an irreducible $(\mathfrak{g},K)$-module. What we are going to do is to attach certain invariants to the representations in $\widehat{G}_{adm}$.\\ 

 \subsection{Infinitesimal Characters} The most basic invariant is the \textit{infinitesimal character} of a representation.  The center $Z(\mathfrak{g})$ of $U(\mathfrak{g})$ can be identified with the $W$-invariant polynomials on $\mathfrak{h}$ via the Harish-Chandra homomorphism $\zeta : Z(\mathfrak{g})  \to U(\mathfrak{h})^W$. In this way, we have a map $infchar: \widehat{G}_{adm} \to \mathfrak{h}^{\ast} / W$, and the infinitesimal character of $\pi\in \widehat{G}_{adm}$ is identified with a weight $\lambda\in \mathfrak{h}^{\ast}$. For $\lambda\in  \mathfrak{h}^{\ast}/W,$ we denote by\begin{center} $\widehat{G}_{adm,\lambda}=\{\pi\in\widehat{G}_{adm} |  infchar(\pi) =\lambda \}$ \end{center} and refer to the representations in $\widehat{G}_{adm,\lambda}$ as the irreducible admissible representations with infinitesimal character $\lambda$. Similarly, let $\mathcal{HC}(\mathfrak{g}, K)_{\lambda}$ denote the set of Harish-Chandra modules with infinitesimal character $\lambda$.
 
\subsection{Primitive Ideals}
Many invariants to be considered are actually invariants attached to the primitive ideals in $U(\mathfrak{g})$, though there are some invariants attached directly to an irreducible Harish-Chandra module. Thus let us first define 

\begin{definition}
Let $V$ be an irreducible $U(\mathfrak{g})$-module. The annihilator of $V$ in $U(\mathfrak{g})$ is\begin{center} Ann$(V) :=\{ X\in U(\mathfrak{g})  |Xv=0,\forall v\in V \}$, \end{center}which is a two-sided ideal in $U(\mathfrak{g})$. It is called the \textit{primitive ideal} in $U(\mathfrak{g})$ attached to $V$. 
\end{definition}
If two $U(\mathfrak{g})$-modules  have the same primitive ideals, then their infinitesimal characters are the same, and hence it makes sense to talk about the primitive ideals with infinitesimal character  $\lambda$. We set
Prim$(\mathfrak{g})_{\lambda}$ to be the set of primitive ideals in $U(\mathfrak{g})$  with infinitesimal character $\lambda$. For any $\pi\in \widehat{G}_{adm, \lambda}$, let $V_{\pi}$ be the corresponding Harish-Chandra module and let $I_{\pi}:=$Ann$(V_{\pi} )$, 
and hence we have a map  $ \widehat{G}_{adm, \lambda}\to \text{Prim}(\mathfrak{g})_{\lambda} $ sending $\pi$ to $I_{\pi}$. This map is several-to-one in general. 

\subsection{Associated  Variety and Gelfand-Kirillov Dimension} \label{AVGKdim}
    
    Given a finitely generated $\mathfrak{g}$-module $V$.
Let $U_n(\mathfrak{g})\subseteq U(\mathfrak{g})$ be the subspace of $U(\mathfrak{g})$ generated by the monomial of the form $X_1\cdots X_m$ with $m\le n$ and $X_i\in \mathfrak{g}$. There is a good filtration (see Section 4 in \cite{BB}) of $V$ compatible with the graded action of $U(\mathfrak{g})$, i.e.  $0=V_0\subset V_1\subset V_2\subset \cdots \subset V$ and $U_n(\mathfrak{g}) V_i\subseteq  V_{i+n}$ for all $i, n$. Then gr$(V)=\bigoplus\limits_{n>0} V_n/V_{n-1}$ is a finitely generated module for the associated graded algebra of $U(\mathfrak{g})$, namely the symmetric algebra $S(\mathfrak{g})$ by Poincar\'e-Birkhoff-Witt theorem. So one can define the \textit{associated variety} of $V$, denoted AV$(V)$, to be the support of the $S(\mathfrak{g})$-module gr$(V)$ in $\mathfrak{g}^{\ast}$.

       Moreover, let $\varphi_V(n)=\sum\limits_{q\le n} \dim_{\mathbb{C}} V_q$, which is finite since $V$ is finitely generated. By a theorem of Hilbert and Serre, there is a polynomial $\overline{\varphi}_V(n)$, of degree at most $\dim \mathfrak{g}$, such that $\varphi_V(n)=\overline{\varphi}_V(n)$ for large $n$. (The proof can be found in \cite{V2}). Therefore, the integer $\deg ( \overline{\varphi}_V)$ is defined to be the \textit{Gelfand-Kirillov dimension} of $V$, denoted by GK$\dim (V).$

         An important lemma is stated below.
\begin{lemma} \label{AV}
AV$(V\otimes F)=$AV$(V)$  and GK$\dim (V\otimes F)=$GK$\dim(V)$ for any finite-dimensional $\mathfrak{g}$-module $F$.
\begin{proof}
Choose a good filtration $\{V_i\}$ on $V$, then we obtain a good filtration $\{V_i\otimes F\}$ on $V\otimes F$. With these filtrations, gr$(V\otimes F)$ as a $S(\mathfrak{g})$-module is a sum of copies of gr$(V)$. Hence the lemma follows. 
\end{proof}
\end{lemma}
Now suppose $\pi \in \widehat{G}_{adm}$ and 
 $I_{\pi}$ is the primitive ideal attached to $\pi$, which can be regarded as a left $U(\mathfrak{g})$-module, and hence we define AV$(I_{\pi})$ and GK$\dim(I_{\pi})$, GKdim$(\pi)$ in usual sense, whereas AV$(\pi)$ will be defined upon a $K$-invariant filtration, and we won't talk about this until Section 4. By Kostant's theory of harmonics, AV$(I_{\pi})$ consists of nilpotent elements in $\mathfrak{g}^{\ast}$, and hence is a union of finite number of closures of nilpotent coadjoint orbits. In fact,  it is a single orbit. Let us  record some remarkable facts as follows.
 
\begin{theorem} \label{uniquecx}
(1) (Borho, Brylinski, Joseph, see \cite{BB} and \cite{J}) There exists a unique (complex) nilpotent coadjoint  orbit $\mathcal{O}$ such that AV$(I_{\pi})=\overline{\mathcal{O}}$.\\
(2) (See \cite{BV}, Theorem 4.1) 2GK$\dim(\pi)=$GK$\dim(I_{\pi})=\dim _{\mathbb{C}} \overline{\mathcal{O}}$, where $\overline{\mathcal{O}}=$ AV$(I_{\pi})$ is obtained from (1).
\end{theorem}

\begin{remark}
AV$(I_{\pi})$, which it the closure of a complex nilpotent orbit, is called the \textit{complex associated variety} of $\pi$, while AV$(\pi)$ (which will be discussed in detail later in Section 4) is the \textit{real associated variety} of $\pi$. 
\end{remark}

\subsection{$\tau$-invariant}
Given $I \in$Prim$(\mathfrak{g})_{\lambda}$. Put $\triangle(\lambda)=\{\alpha\in\triangle | \la\lambda, \alpha^{\vee}\ra\in \mathbb{Z}   \}$, the integral root system for $\lambda$, and let $W_{\lambda}$ denote the Weyl group for $\triangle (\lambda)$. Choose $\triangle^+ (\lambda) \subseteq  \triangle(\lambda)$, a positive system making $\lambda$ dominant.  Write $\prod (\lambda) \subseteq \triangle^+(\lambda)$ for the set of simple roots. There is the Borho-Jantzen-Duflo $\tau$-invariant attached to $I$, which is a subset of $\prod (\lambda)$ (see \cite{V3}), denoted $\tau (I)$. 

  Since $G_{\mathbb{C}}$ is simply connected, we have an alternative definition for  $\tau$-invariant (see \cite{IC2}, Section 2).  Let $\pi\in \mathcal{HC}(\mathfrak{g},K )_{\lambda}$ and $F_{\gamma}$ be the finite-dimensional representation of $G$ with extreme weight $\gamma$. Also let  $\triangle(F_{\gamma})$ denote the set of all weights of $F_{\gamma}$. Consider the Zuckerman translation functor $\psi _{\lambda} ^{\lambda +\gamma} (\pi)=P _{\lambda+\gamma} ({\pi \otimes F_{\gamma}})$, where $P_{\lambda +\gamma}$ by definition is the projection on the representations with infinitesimal character $\lambda +\gamma$, and hence  $\psi _{\lambda} ^{\lambda +\gamma} (\pi)$ is a functor that  projects $\pi \otimes F_{\gamma} $ on representations with infinitesimal character $\lambda +\gamma$.  Let $\alpha \in \prod (\lambda)$, and let $\lambda _{\alpha}$  be singular with respect to $\alpha$ such that $\lambda-\lambda _{\alpha}$ is a sum of roots. Define $\psi _{\alpha}(\pi) := \psi _{\lambda}^{\lambda_{\alpha}} (\pi)$ be the translation functor of $\pi$ to the $\alpha$-wall. Then we define \begin{center} $\tau(\pi) =\{\alpha \in \prod (\lambda) | \thinspace \psi_{\alpha} (\pi) =0 \}$. \end{center}  It turns out that   $\tau$-invariant is a measure of size of $\pi$: the bigger the $\tau$-invariant, the smaller the representation.  
\begin{definition}
We say that $\pi$ has maximal $\tau$-invariant if $\tau(\pi)=\prod (\lambda)$,  or equivalently, $\psi _{\alpha} (\pi)=0$ for all $\alpha\in \prod (\lambda)$.
\end{definition}
\begin{lemma}  \label{findimmaxtau}
Let $F$ be a finite dimensional representation. Then $\psi_{\alpha} (F)=0$ for every root $\alpha$ and hence $F$ has maximal $\tau$-invariant.
\begin{proof}
Note that the infinitesimal character of every finite dimensional representation is regular. 

    Assume the setting in the Lemma. We have $\psi_{\alpha}(F)=P_{\lambda'}(F\otimes F')=0$, where $\lambda'$ is singular for $\alpha$ and $F'$ is a finite dimensional representation, since $F\otimes F'$ is a virtual finite dimensional representation and each constituent has regular infinitesimal character.  
\end{proof}
\end{lemma}

\begin{definition}  \label{small}
We call a representation \textit{small} if it has maximal $\tau$-invariant. 
\end{definition}
The Gelfand-Kirillov dimension of an irreducible representation is a measure of the growth of $K$-types. Here is the proposition connecting these two measures.

\begin{proposition} (\cite{V1}, Corollary 4.5) Let $\pi\in\widehat{ G}_{adm, \lambda}$. If $I_{\pi}$ has max $\tau$-invariant, then \begin{center} GK$\dim (\pi)=|\triangle ^+|-|\triangle ^+ (\lambda)|.$\end{center}
\end{proposition}

\subsection{Weyl Group Representations}
There are some details of Weyl group representations that can be found in various places, for instance, \cite{Ca}, \cite{Lus}, \cite{McG},  and \cite{T1}. We recall some of the useful facts as follows.

  In \cite{J2}, Joseph has attached to $I\in$Prim$(\mathfrak{g}) _{\lambda}$ a representation $\sigma_{I}\in \widehat{W}_{\lambda} $. In fact, the map from $I\in$ Prim$(\mathfrak{g}) _{\lambda}$ to  $\sigma_{I}$ is surjective onto the set of special representations of $W_{\lambda}$ (see \cite{Ca} for definition of a special Weyl group representation). 

   On the other hand, Springer provides a method for producing a representation of $W$ from a nilpotent orbit $\mathcal{O}$, which is the well-known Springer correspondence. We write sp$(\mathcal{O})$ for the irreducible representation of $W$ attached to $\mathcal{O}$. There is an algorithm to calculate  $sp(\mathcal{O})$ if given $\mathcal{O}$ by use of symbols (see \cite{McG}). Note that the map $\mathcal{O}\to$ sp$(\mathcal{O})$ is injective, but not surjective usually.
   
    Let $W'$ be any subgroup of W generated by reflections. There is an operation called  t\textit{runcated induction} $j_{W'} ^W$ (see \cite{Ca}, \cite{Lus}), taking irreducible representations of $W'$ to those of $W$.   It is a fact that  $j_{W'} ^W$ is an injective map.

The following proposition summarizes and connects all concepts stated above.
\begin{proposition}   
Let  $\pi\in \widehat{G}_{adm, \lambda}$, $I=I_{\pi}$, $W_{\lambda}$ be the integral Weyl group for $\lambda$. Then
there is a unique nilpotent orbit $\mathcal{O}$ such that $\sigma _I =$sp$(\mathcal{O)}$. Furthermore, this $\mathcal{O}$ is dense in AV$(I)$, that is, AV$(I)=\overline{\mathcal{O}}$. Thus, we  have a commutative diagram: 
\begin{center}
\begin{tikzpicture}
\tikzset{node distance=2cm, auto}
  \node (P) {$\mathcal{O}$};
  \node (B) [right of=P] {sp$(\mathcal{O})$};
  \node (A) [below of=P] {$I$};
  \node (C) [below of=B] {$\sigma_I$};
  \draw[->] (A) to node {AV} (P);
  \draw[->] (P) to node [swap] {sp} (B);
  \draw[->] (A) to node [swap] {$\sigma$} (C);
  \draw[->] (C) to node {$j$} (B);
\end{tikzpicture} 
\end{center}
(The left vertical arrow in the diagram means $AV(I) = \overline{\mathcal{O}}$.)
\begin{proof}
The proof can be found in Section 2.5 and 2.6 in \cite{T1}.
\end{proof}
\end{proposition}

\section{ Some Small Representations of $\widetilde{G}$} \label{sectionsmall}

In this section we assume that $G$ is a real form of a simply connected, semisimple complex Lie group, and $\widetilde{G}$ is the nonlinear double cover of $G$.  First, we identify the kernel of the covering map $p: \widetilde{G}\to G$ with $\pm 1$ and write $\widetilde{H}$ for the inverse image in $\widetilde{G}$ of a subgroup $H$ of $G$. We define 

\begin{definition}
A representation $\widetilde{\pi}$ of $\widetilde{H}$ is called genuine if $\widetilde{\pi}(-1)=-I$. If $\widetilde{\pi}$ is irreducible, then $\widetilde{\pi}$ is genuine if and only if $\widetilde{\pi}$ does not factor through $H$.
\end{definition}

We focus on the genuine representations with a particular infinitesimal character $\lambda$.  As mentioned in the introduction, such $\lambda$ is chosen such that the irreducible quotients of the pseudospherical principal 
series representations defined  in \cite{ABPTV} are small in the sense that they have maximal $\tau$-invariants. Indeed, $\lambda$  can be defined by the following formula, 
\begin{equation}
\lambda =\frac{1}{2} \left (  \sum _{\alpha \text{ not metaplectic }}  \alpha \right ) + \frac{1}{4} \left (  \sum _{\alpha \text{ metaplectic }}  \alpha \right ),
\end{equation}
where the notion of a root being metaplectic is defined in \cite{ABPTV}.  The infinitesimal characters for all types are 
  listed in Table 1.  
  
  We are interested in a special category of representations with certain properties, defined as follows.\\ 
 Denote\begin{center} $\prod _{\lambda} ^s (\widetilde{G})=\{ \widetilde{ \pi} \thinspace |   \thinspace \widetilde{\pi} \in \widehat{\widetilde{G}}_{adm,\lambda} \text{, $\widetilde{\pi}$  is  genuine and has maximal $\tau$-invariant} \}$ \end{center}
The following is the key Lemma.

\begin{lemma} \label{cxnilp}
There is a unique complex nilpotent orbit $\mathcal{O}$ such that AV$(I_{\widetilde{\pi}})=\overline{\mathcal{O}}$ for every $\widetilde{\pi}\in\prod _{\lambda} ^s (\widetilde{G})$. This $\mathcal{O}$ can be computed explicitly (case by case) and 
 $\dim ( \mathcal{O})=2$GK$\dim(\widetilde{\pi})=2(|\bigtriangleup ^+|-|\bigtriangleup ^+(\lambda)|)$, where $\triangle$ and $\triangle(\lambda)$ are  the root system and integral root system, respectively. 

\begin{proof}
Let $\pi\in\prod _{\lambda} ^s (\widetilde{G})$. Since $\widetilde{\pi}$ has maximal $\tau$-invariant, $\sigma_{I_{\widetilde{\pi}}} =sgn_{W_{\lambda}}$, the sign representation of the integral Weyl group for $\lambda$. Then the truncated induction takes $sgn_{W_{\lambda}}$ to a special representation of $W$, denoted $j(sgn)=j_{W_{\lambda}} ^{W}  (sgn)$, since $sgn_{W_{\lambda}}$ is a special representation of $W_{\lambda}$. Hence $j (sgn)$ defines a nilpotent orbit $\mathcal{O}$ of $\mathfrak{g}$ through the Springer correspondence, i.e. $sp(\mathcal{O})=j(sgn)$, and this $\mathcal{O}$ is dense in the associated variety of $I_{\widetilde{\pi}}$, which means AV$(I_{\widetilde{\pi}})=\overline{\mathcal{O}}$.  The uniqueness of this $\mathcal{O}$ follows from either Theorem \ref{uniquecx} (1) or the injectivity of the Springer correspondence.  

    From \cite{V1}, for a representation at infinitesimal character $\lambda$ with maximal $\tau$-invariant, GK$\dim (\widetilde{\pi})=|\bigtriangleup ^+|-|\bigtriangleup ^+(\lambda)|$, and hence $\dim \mathcal{O} = 2(|\bigtriangleup ^+|-|\bigtriangleup ^+(\lambda)|)$ by Theorem \ref{uniquecx} (2). For exceptional groups, there is a unique complex nilpotent orbit of this dimension (see \cite{CM}), so it is exactly the one that we are looking for.  
    
    For classical types,  there is an algorithm to calculate $j(sgn)$ and the corresponding $\mathcal{O}$ explicitly  via the Springer correspondence (see \cite{CM}). The  parametrization sets of nilpotent orbits  are partitions of $n$ for type $A_{n-1}$,  and are partitions of $2n (\thinspace 2n+1,\text{ resp.})$  which even (odd, resp.) parts occur with even  multiplicity  for type $B_n$ and $D_n$ ($C_n$, resp.) (see \cite{Ca} and \cite{CM}). All of the nilpotent orbits and the corresponding Weyl group representations are listed in Table 1.
\end{proof}

\end{lemma}
\begin{table}[H] \label{table1}
\caption{}
\begin{tabular}{ccccccc} Type & $n$ & $\lambda$ &$\triangle(\lambda)$  &$\dim \mathcal{O}$&$ \mathcal{O}$ & $j(sgn_{W_{\lambda}})$ \\\hline \multirow{2}{*}{ $A_{n-1} $}  &$2m$ &$\rho/2$ & $A_{m-1}\times A_{m-1}$ & $\frac{n^2}{2}$&$[2^m]$ & $[2^m]$  \\  \cline{2-2}\cline{3-3}   \cline{4-4}\cline{5-5} \cline{6-6} \cline{7-7}
  & $2m+1$ & $ \rho/2$ & $A_{m-1}\times A_{m-1}$ & $\frac{n^2-1}{2}$ & $[2^m\thinspace 1]$ &$[2^m\thinspace 1]$  \\\hline 
   \multirow{2}{*}{ $B_n $}  &$2m$ &$\rho(C_n)/2$ & $B_m\times B_m$ & $n^2$&$[2^n\thinspace 1]$ & $(\phi; [2^m])$  \\  \cline{2-2}\cline{3-3}   \cline{4-4}\cline{5-5} \cline{6-6} \cline{7-7}
  & $2m+1$ & $ \rho(C_n)/2$  & $B_{m+1}\times B_m$  & $n^2-1$ & $[2^{n-1}\thinspace 1^3]$ &$(\phi; [2^m\thinspace 1])$   \\\hline $C_n$ && $ \rho(B_n)$& $D_n$&$2n$& $[2\thinspace 1^{2n-2}]$& $([1^n];\phi)$\\\hline
  \multirow{2}{*}{ $D_n $}  &$2m$ & $\rho/2$& $D_m\times D_m$  & $n^2$&$[3\thinspace 2^{n-2}\thinspace 1]$ & $\{\phi; [2^m]\}$  \\  \cline{2-2}\cline{3-3}   \cline{4-4}\cline{5-5} \cline{6-6} \cline{7-7}
  & $2m+1$& $\rho/2$ & $D_{m+1}\times D_m$  & $n^2-1$ & $[3\thinspace 2^{n-3}\thinspace 1^3]$ &$\{\phi; [2^m\thinspace 1]\}$   \\\hline   $E_6$& & $\rho/2$&$A_1\times A_5$&40&$3A_1$& $\phi_{15,16}$\\\hline  $E_7$&&$\rho/2$&$A_7$&70&$4A_1$&$\phi_{15,28}$\\\hline $E_8$ & &$\rho/2$&$D_8$&128&$4A_1$&$\phi_{50,56}$ \\\hline  $F_4$&&$\rho ^{\vee}$&$B_4$&16&$A_1$&${\phi _{2,16} }''$ \\\hline $G_2$ &&$\rho/2$&$A_1\times A_1$&8&$\widetilde{A_1}$&$\phi_{2,2}$ 

 \end{tabular}
\end{table}

Because of Lemma \ref{cxnilp}, let $\mathcal{O}$ denote the complex nilpotent orbit such that AV$(I_{\widetilde{\pi}})=\overline{\mathcal{O}}$ for $\widetilde{\pi} \in \prod_{\lambda} ^s (\tu{G})$, and define
\begin{center}$\prod _{\lambda} ^{\mathcal{O}} (\widetilde{G})= \{\widetilde{\pi} \thinspace |   \thinspace\widetilde{ \pi }\in \widehat{\widetilde{G}}_{adm,\lambda}, \widetilde{\pi}\text{ is genuine and    AV$(I_{\widetilde{\pi}})=\overline{\mathcal{O}}$}\}$.
\end{center}
Then here is the main theorem of this section.
\begin{theorem} \label{main1} $\prod _{\lambda} ^s (\widetilde{G})= \prod _{\lambda} ^{\mathcal{O}} (\widetilde{G})$.
\begin{proof}
It is clear that $\prod _{\lambda} ^s (\widetilde{G}) \subseteq \prod _{\lambda} ^{\mathcal{O}} (\widetilde{G})$ due to Lemma \ref{cxnilp}. Conversely, given a representation $\widetilde{\pi}\in  \prod _{\lambda} ^{\mathcal{O}} (\widetilde{G}) $, we need to show that $\widetilde{\pi}$ has maximal $\tau$-invariant, that is, to show that $\sigma _{I_{\widetilde{\pi}}} = sgn _{W_{\lambda}}$. This  simply follows from the injectivity of the truncated induction. 
\end{proof}
\end{theorem}

\section{Real Associated Variety}
In the previous section, given $\pi\in \widehat{G }_{adm}$, we defined its complex associated variety AV$(I_{\pi})$. Now we want to attach nilpotent orbits directly to $\pi$.  Notice that these notions are quite general and they can be defined for linear and nonlinear groups.

   Suppose $(\pi,V)$ is the given finitely-generated $(\mathfrak{g} , K)$-module. As in Section \ref{AVGKdim}, suppose  $0=V_0\subset V_1\subset V_2\subset \cdots \subset V$ is a good filtration, and furthermore suppose this is a  $K_{\mathbb{C}}$-invariant  filtration since $V$ is also a $K_{\mathbb{C}}$-module. Hence the associated variety of $(\pi, V)$, denoted AV$(\pi)=$AV$(V)$, is a closed subvariety of $(\mathfrak{g}/\mathfrak{k})^{\ast}$. Since $V$ is also a $K_{\mathbb{C}}$-module, AV$(\pi)$ is actually a $K_{\mathbb{C}}$-invariant subset of $(\mathfrak{g}/\mathfrak{k})^{\ast}$. Similarly, AV$(\pi)$ consists of nilpotent elements, say, AV$(\pi)\subseteq \mathcal{N}(\mathfrak{g}/\mathfrak{k})^{\ast} :=\mathcal{N}(\mathfrak{g}^{\ast})\cap (\mathfrak{g}/\mathfrak{k})^{\ast}$, where $\mathcal{N}(\mathfrak{g}^{\ast})$ denotes the nilpotent cone of $\mathfrak{g}^{\ast}$.  By a theorem of  Kostant-Rallis (see Proposition 4 in \cite{KR}), there are finitely many $K$ orbits on $\mathcal{N}(\mathfrak{g}/\mathfrak{k})^{\ast} $, and hence we may write \begin{center} AV$(\pi)=\overline{\mathcal{O}_1^{K_{\mathbb{C}}}}\cup\cdots\cup \overline{\mathcal{O}_j ^{K_{\mathbb{C}}}}$,\end{center} for orbits $\mathcal{O}_i ^{K_{\mathbb{C}}}$ of $K_{\mathbb{C}}$ on $ \mathcal{N}(\mathfrak{g}/\mathfrak{k})^{\ast}.$
 
     The next result of Vogan relates AV$(I_{\pi})$ and AV$(\pi)$.   
\begin{theorem}   (see \cite{V4}, Section 5) Suppose $\pi\in \widehat{G}_{adm}$. Write     \begin{center} AV$(\pi)=\overline{\mathcal{O}_1^{K_{\mathbb{C}}} }\cup\cdots\cup \overline{\mathcal{O}_j ^{K_{\mathbb{C}}}}$, and AV$(I_{\pi})=\overline{\mathcal{O}}$. \end{center}  Then each $\mathcal{O}_i^{K_{\mathbb{C}}}$ is a Lagrangian submanifold of the canonical symplectic structure of $\mathcal{O}$. In particular, for each $i$, we have  \begin{center}  $G\cdot\mathcal{O}_i^{K_{\mathbb{C}}}=\mathcal{O}$ and GK$\dim (\pi)=\dim (\mathcal{O}_i ^{K_{\mathbb{C}}}  )$.   \end{center}
\end{theorem}

Next we introduce the Sekiguchi correspondence (see \cite{CM}, Chapter 9, for example). 

\begin{theorem} (Sekiguchi)
There is a natural one-to-one correspondence between nilpotent $G$-orbits in $\mathfrak{g}_{\mathbb{R}}$ and nilpotent $K_{\mathbb{C}}$-orbits in $(\mathfrak{g}/\mathfrak{k})$.
\end{theorem}

Thus, by the Sekiguchi correspondence, AV$(\pi)$ can be viewed as $\overline{\mathcal{O}_1}\cup\cdots \overline{\mathcal{O}_j}$, where each $\mathcal{O}_i$ is a $G$-orbit in $\mathfrak{g}_{\mathbb{R}}$ corresponding to $\mathcal{O}_i ^{K_{\mathbb{C}}}$ via the Sekiguchi correspondence, and hence AV$(\pi)$ is called the\textit{ real associated variety} of $\pi$. Moreover, if  AV$(I_{\pi})=\overline{\mathcal{O}}$, then we have   $G_{\mathbb{C}}\cdot \mathcal{O}_i=\mathcal{O}$, and hence we say that each $\mathcal{O}_i$ is a real form of $\mathcal{O}$. Equivalently, we say $\{\mathcal{O}_i\}_{i=1}  ^l$ is the set of real forms of $\mathcal{O}$ if $\mathcal{O} \cap \mathfrak{g}_{\mathbb{R}}= \mathcal{ O}_1 \cup \cdots \cup \mathcal{O}_l.$

\subsection{Real associated variety of representations in $\prod_{\lambda} ^s (\tu{G})$}
Resuming the setting of $G$ and $\widetilde{G}$ in Section  \ref{sectionsmall}, recall that we defined a set of representations $\prod ^s _{\lambda} (\widetilde{G})$,   and the complex associated variety of each representation in this set is the closure of a particular $\mathcal{O}$ (see Table 1). In Table 2, we list all real groups $G$ such 
that $\mathcal{O}  \cap \mathfrak{g}_{\mathbb{R}}$ is nonempty,  as well as the number of real forms of $\mathcal{O}$, denoted $\# \mathcal{O}_i$,  with respect to each $G$ (see \cite{CM} for the parametrization of real nilpotent orbits).

\begin{remark} \label{realorbit}
It can be observed from Table 2 that there are not many real groups which have nonempty intersection with $\mathcal{O}$. More precisely, if $G$ is not listed in Table 2, then $\mathcal{O}\cap\mathfrak{g}_{\mathbb{R}}=\phi$. 
\end{remark}

\begin{table}[H] \label{table2}
\caption{}
\begin{center}
\begin{tabular}{|c|c|c|c|c|c|c|c|}\hline
Type &\multicolumn{4}{c|}{$A_{n-1} \quad (\mf{g}=\mf{sl}_n$)} &  \multicolumn{3}{c|}{$B_{n} \quad(\mf{g}=\mf{so}_{2n+1}$)} \\
\hline
$n$&$2m$& $2m+1$&$2m$& $2m+1$ & $2m$&$2m+1$&$2m+1$ \\\hline $G$&  \multicolumn{2}{c|}{$SL(n,\R)$}& $SU(m,m)$&$SU(m+1,m)$ &    \multicolumn{2}{c|}{$Spin(n+1,n)$} &$Spin(n+2,n-1)$ \\\hline $\# \Or _i$&2&1&$m+1$& $m+1$&2&1&1\\\hline
\end{tabular} 
\end{center}
\end{table}
\begin{table}[H] 
%\caption{$\mathcal{O}$ and its corresponding Weyl group representation (using parameterizations in \cite{Ca} and \cite{CM})}
\begin{center}
\begin{tabular}{|c|c|c|c|c|c|c|c|}\hline
Type &\multicolumn{2}{c|}{$C_{n} \quad (\mf{g}=\mf{sp}_{2n}$)} &  \multicolumn{5}{c|}{$D_{n}\quad (\mf{g}=\mf{so}_{2n}$)} \\
\hline
$n$&& &$2m$& $2m+1$ & $2m$&$2m+1$&$2m+1$ \\\hline $G$& $Sp(2n,\R)$&$Sp(2p,2q)$& \multicolumn{2}{c|}{$Spin(n,n)$} &    \multicolumn{2}{c|}{$Spin(n+1,n-1)$} &$Spin(n+2,n-2)$ \\\hline $\# \Or _i$&2&1&4&2&1&2&1\\\hline
\end{tabular} 
\end{center}
\end{table}
\begin{table}[H] 
%\caption{$\mathcal{O}$ and its corresponding Weyl group representation (using parameterizations in \cite{Ca} and \cite{CM})}
\begin{center}
\begin{tabular}{|c|c|c|c|c|c|c|}\hline
Type &\multicolumn{2}{c|}{$E_6$} &  $E_7$&$E_8$&$F_4$&$G_2$ \\
\hline
$G$& $E_6(A_1\times A_5)$&$E_6(C_4)$& $E_7(A_7)$&$E_8(D_8)$&$F_4(B_4)$& $G_2(A_1\times A_1)$ \\\hline $\# \Or _i$ &2&1&2&1&1&1\\\hline
\end{tabular} 
\end{center}
\end{table}

We have the following proposition saying that we can attach to each small representation defined in Section \ref{sectionsmall} a single real nilpotent orbit.

\begin{proposition} \label{singlereal}
We resume the  setting  and notations in Section \ref{sectionsmall}. Suppose $G_{\mathbb{C}}$ is a simply connected, semisimple complex Lie group, $G$ is a real form of $G_{\mathbb{C}}$ and $\widetilde{G}$ is the nonlinear double cover of $G$. For each $\widetilde{\pi}\in \prod _{\lambda} ^s (\widetilde{G}) =\prod _{\lambda} ^{\mathcal{O}} (\widetilde{G})$, there is a unique real nilpotent orbit $\mathcal{O}_{\widetilde{\pi}}$  such that AV$(\widetilde{\pi})=\overline{\mathcal{O}_{\widetilde{\pi}}}$. This $\mathcal{O}_{\widetilde{\pi}}$ is one of the real forms of $\mathcal{O}$.
\begin{proof} By a result of Vogan (see Theorem 1.3 in  \cite{V4}), if $\mathcal{O}_l$ is a real orbit of maximal dimension in AV$(\widetilde{\pi})$, and the complement of $\mathcal{O}_l$ has codimension at least two in $\overline{\mathcal{O}_i}$, then AV$(\widetilde{\pi})=\overline{\mathcal{O}_l}$. Since $\dim _{\mathbb{R}} \mathcal{O}_i =\dim _{\mathbb{C}} \mathcal{O}$ for each real form $\mathcal{O}_i$ of $\mathcal{O}$, we just need to pick a complex nilpotent orbit $\mathcal{O}'$, which is one step down smaller than $\mathcal{O}$, and see if the difference of $\dim \mathcal{O}$ and $\dim \mathcal{O}'$ is at least 2. This can be checked case by case.
\end{proof}
\end{proposition}

\section{Genuine Central Characters}

\subsection{Regular Characters}\label{reg}
The following material can be found in several places, for example, \cite{AT},  \cite{RT}, and \cite{V5}. Again, $G$ is a real form of a  simply connected, semisimple complex Lie group and $\widetilde{G}$ is the nonlinear double cover of $G$.  Let $\pi\in \widehat{G}_{adm, \lambda}$, where $\lambda$ is a regular infinitesimal character. Then $\pi$ can be specified by a parameter, which is called a $\lambda$-regular character,  $\gamma=(H,\Gamma,\overline{ \gamma})$, where $H$ is a $\theta$-stable Cartan subgroup of $G$, $\Gamma$ is a character of $H$, and $\overline{\gamma}$ is an element in $\mathfrak{h}^{\ast}$ which defines the same infinitesimal character as $\lambda$, and there are  certain compatibility conditions between $ \overline{\gamma}$  and $\Gamma$ (see Definition 5.3 in \cite{AT}). Write $H=TA$, where $T=H^{\theta}$ and $A$ is the identity component of $\{ h\in H | \theta (h)=h^{-1}\}$. Let $M=Cent_G(A)$. The conditions on $\gamma$ imply that there is a unique relative discrete series representation of $M$, denoted by $\sigma_M$,
with Harish-Chandra parameter $\overline\gamma$, whose lowest $M\cap K$-type has $\Gamma$ as a highest weight.
Then we define a parabolic subgroup $P=MN$ such that $\pi=J(\gamma)$, the unique irreducible quotient of a standard representation $I(\gamma)= Ind _P ^G (\sigma _M \otimes 1)$, which is parametrized by $\gamma$ from a $K$-conjugacy class of regular characters for $\lambda$. 

Recall that (see \cite{A1}, for instance) when $\lambda$ is a regular infinitesimal character, $\mathcal{HC}(\mathfrak{g}, K )_{\lambda}$ is parametrized by the set $\mathcal{P}_{\lambda}$ of $K$-conjugacy classes of $\lambda$-regular characters. Furthermore, the following two sets are bases of the Grothendieck group: \begin{center} 
$\{[ J(\gamma)   ] \} _{\gamma\in\mathcal{P}_{\lambda} }$ and $\{[ I(\gamma)   ] \} _{\gamma\in\mathcal{P}_{\lambda} }$. \end{center}

We have the following definition.
\begin{definition} \label{KLVpoly}
Define the change of basis matrix \begin{center} $[J (\delta)]  = \sum \limits _{\gamma \in \mathcal{P}_{\lambda}  }  M(\gamma , \delta)  [I(\gamma)]$  \end{center}
and the inverse matrix \begin{center} $[I   (\delta)]  = \sum \limits _{\gamma \in \mathcal{P}_{\lambda}  }  m(\gamma , \delta)  [J(\gamma)]$.  \end{center}

Here $M(\gamma , \delta)$ and $m(\gamma , \delta)$ are integers and $M(\gamma , \delta) $ are computed by the Kazhdan-Lusztig-Vogan algorithm when $G$ is linear.
\end{definition}

%In particular, consider $\mathbb{C}$, the trivial representation, and write its standard module as $I(\gamma_0)$ with parameter $\gamma_0$. Then the coefficients $M(\gamma, \gamma_0)$ are $\pm 1$. 

%\begin{lemma} \label{charflaC}  (\cite{V5})
%There is an identity in the Grothendieck group \begin{center}  $\mathbb{C}  =  \sum \limits _{\gamma}  (-1)^{l(\gamma _0) -l(\gamma) } I(\gamma),    $\end{center}
%where $\gamma=(H, \Gamma, \overline{\gamma})$ runs over holomorphic characters $\Gamma$ on $H$ (see \cite{A1} for definition), and $l(\gamma)$ is the length function in \cite{V5}.
%\end{lemma}

    The above notions can also be defined for nonlinear groups. More specifically, let $\lambda$ be the regular infinitesimal character  defined in Section \ref{sectionsmall}. In this case, suppose that $\widetilde{\pi}$ is an irreducible genuine representation from $\widehat{\widetilde{G} }_{adm,\lambda}$. Then $\widetilde{\pi}$ is parametrized by a genuine $\lambda$-regular character $\gamma=(\widetilde{H}, \Gamma,\overline{ \gamma})$, where  $\Gamma$ is an irreducible genuine representation of $\widetilde{H} = p^{-1} (H)$. Note that in this case $\Gamma$ can be replaced by a character of $Z(\widetilde{H})$, a central character of $\widetilde{H}$, because of the following proposition (see \cite{ABPTV}).
    
\begin{notation}    
Let $G'$ be a subgroup of $\tu{G}$. We write $Z(G')$ for the center of $G'$ and $\prod _g (G')$ for equivalence classes of irreducible genuine representations of $G'$.
\end{notation}    
\begin{proposition}
Let $\gamma=(\widetilde{H}, \Gamma,\overline{ \gamma})$ be a genuine regular character defined as above.  Let $n=|\widetilde{H}/ Z(\widetilde{H})| ^{\frac{1}{2}}$.  For every 
$\chi \in \prod _g (Z(\widetilde{H}))$ there is a unique representation $\Gamma=\Gamma(\chi)\in \prod _g (\widetilde{H})$ for which $\Gamma | _{Z(\widetilde{H})}$ is a multiple of $\chi$. The map $\chi\to\Gamma(\chi)$ is a bijection between $\prod _g (Z(\widetilde{H}))$ and $\prod _g (\widetilde{H})$. The dimension of $\Gamma(\chi)$ is $n$, and $Ind_{Z(\widetilde{H})} ^{\widetilde{H}} (\chi)  =n\Gamma$.
\begin{proof}

The proof can be found in Proposition 2.2 in \cite{ABPTV}.
\end{proof}
\end{proposition}

When $G$ is simply laced, a genuine representation of $\widetilde{H}$ is determined by the  infinitesimal character and its restriction to $Z(\widetilde{G})$. We record the properties as follows (see \cite{AHe} and  \cite{AT}).

\begin{proposition} \label{center} 
All the setting is as before, and also suppose $G$ is simply laced, $H$ is a Cartan subgroup of $G$, and $H^0$ is the identity component of $H$. Then\\
(1) (\cite{AHe}, Proposition 4.7) $Z(\widetilde{H})=Z(\widetilde{G})\widetilde{H^0}$. In particular, a genuine character of $Z(\widetilde{H})$ is determined by its restriction to $Z(\widetilde{G})$ and its differential;\\
(2) (\cite{AT}, Proposition 5.5) A genuine regular character $\gamma= (\widetilde{H}, \Gamma,\overline{ \gamma})$ of $\widetilde{G}$ is determined by $\overline{\gamma}$ and the restriction of $\Gamma$ to $Z(\widetilde{G})$, and so is $\widetilde{\pi}=J(\gamma)$.
\end{proposition}

The second part of this proposition is basically a corollary of the first part. Consequently $\widetilde{G}$ typically has few genuine irreducible representations, denoted $\prod _g (\widetilde{G})$. 

The following corollary immediately follows from Proposition 5.4, which is also mentioned in the example below Proposition 5.5 in \cite{AT}.

\begin{corollary}\label{c-1ps}
 Suppose $\tu{G}$ is split and fix a genuine central character and an infinitesimal character of $\tu{G}$. Then there is precisely one minimal principal series representation with the given infinitesimal character 
and central character. 
\end{corollary}

\subsection{Action of Aut$(G)$ on $\prod _g (\widetilde{G})$}\label{actionaut}

In this section we want to see how an automorphism of $G$ acts on	 $\prod _g(\widetilde{G})$. Let Aut$(G)$ denote the automorphisim group of $G$,  and 
\begin{center}$\text{Int}(G) = \{ \tau  \in \text{Aut}(G) \thinspace | \thinspace \tau =\tau _x \text{ for some } x \in G\}$, where $
\tau_x (g ) = xgx^{-1}$ for $g\in G,$ \\Out$(G)= \text{Aut}(G) / \text{Int}(G).$ \end{center}

\begin{lemma} \label{outeraction}
There is a natural map from Out$(G)$ to Aut$(Z(\widetilde{G}))$, which sends each $\tau \in $Aut$(G)$ to $\widetilde{\tau}\in$   Aut$(Z(\widetilde{G}))$. 
\begin{proof}
The map $\tau\in$Aut$(G) \to\widetilde{ \tau}\in$Aut$(Z(\widetilde{G}))$ is defined as follows. Every $\tau\in$Aut$(G)$ can be lifted to an automorphism $\widetilde{\tau}$ of $\widetilde{G}$. Then by restricting $\widetilde{\tau}$ to $Z(\widetilde{G})$, we get an automorphism of $Z(\widetilde{G})$, which is also denoted by $\widetilde{\tau}$. This map is well-defined since if $\tau\in$Int$(G)$ , say, $\tau =\tau_x$ for some $x\in G$, $\widetilde{\tau} (\widetilde{z}) = \widetilde{x}\widetilde{z}\widetilde{x} ^{-1}=\widetilde{z}\widetilde{x}\widetilde{x}^{-1}=\widetilde{z}$, for $\widetilde{z}\in Z(\widetilde{G})$. %The proof of the second assertion can be found in \cite{AHe}.
\end{proof}
\end{lemma}

Let $\tau\in$Aut$(G)$. Define an action of $\tau$ on $\prod _g(Z(\widetilde{G}))$ as follows. Let $\chi\in\prod _g(Z(\widetilde{G}))$, define $\chi ^{\tau} (z) := \chi (\widetilde{\tau} (z) )$, $z\in Z(\widetilde{G})$.  When $G$ is simply laced, we have an action of Aut$(G)$ on $\prod _g (\widetilde{G})$. Due to Proposition \ref{center} (2), every $\pi = J(\gamma)$, where $\gamma=(\widetilde{H}, \Gamma,\overline{ \gamma})$,  is determined by $\overline{\gamma}$ and $\chi :=\Gamma |_{Z(\widetilde{G})} $. Then we can define $\pi ^{\tau} := J(\gamma ^{\tau}) $, where $\gamma ^{\tau}$  is a regular character determined by $\overline{\gamma}$ and $\chi^{\tau}$.

\begin{remark}
 Suppose $G$ is simply laced. Let $\tau\in$Aut$(G)$ and  $\widetilde{\pi} = J(\gamma)\in \prod_g( \widetilde{G})$, where $\gamma =(\tu H, \Gamma, \overline \gamma)$ with $\chi = \Gamma \mid _{Z(\tu G)}$.  If $\chi \neq \chi ^{\tau}$,  then $\widetilde{\pi}$ and $\widetilde{\pi}^{\tau}$  are inequivalent representations in $\prod _g (\widetilde{G})$.

\end{remark}

\section{ $\prod _{\lambda} ^s (\widetilde{G})$ -- Split Case}\label{splitcase}

In this section, the setting is as in Section \ref{sectionsmall}. Furthermore we assume that $G$ is simply laced and split, and hence $\lambda =\rho/2$ from now on. It follows from Corollary \ref{c-1ps} that there is a unique minimal principal representation of $\widetilde{G}$ if  we fix a genuine central character and infinitesimal character. In \cite{ABPTV}, the unique irreducible quotients of them are called irreducible genuine pseudospherical representations. We will show in this section that we can get more representations in $\prod _{\rho/2} ^s (\widetilde{G})$  by
applying  Cayley transforms to the irreducible genuine pseudospherical representations with infinitesimal character $\rho/2$.

\subsection{Irreducible genuine pseudospherical representations}
Let $A$ denote the split Cartan subgroup of $G$ with Lie algebra $\mathfrak{a}$. 
Let $B$ be a Borel subgroup of $G$ with unipotent radical $N$. We write $B=MA^0N$, where $A^0$ is the identity component of $A$ and $M=Cent_K(A^0)$. The group  $\tu{B}=p^{-1}(B)$ has the semidirect product decomposition: $\tu{B}=\tu{M}A^0N$. 
In \cite{ABPTV},  the pseudospherical principal series representations are defined as follows:
$$
I(\tu{\delta}, \nu ) := Ind ^{\tu{G}} _{\tu{M}A^0N} (\tu{\delta} \otimes e^{\nu} \otimes 1),
$$
where $\tu{\delta}$ is an irreducible genuine pseudospherical representation of $\tu{M}$ (see Definition 4.9 in \cite{ABPTV}). We specify $\nu =\rho/2$. The unique irreducible quotient $J(\tu{\delta}, \rho/2)$ of $I(\tu{\delta}, \rho/2)$ 
are called irreducible genuine pseudospherical  representations with infinitesimal character $\rho/2$. By Proposition 5.6 in \cite{ABPTV}, there is a bijection between the set of irreducible genuine pseudospherical representations of 
$\tu M$ and the set of genuine pseudospherical principal series of $\tu G$ with fixed infinitesimal character. Therefore, we denote 
$$\SH _{\rho/2}=\{ J(\tu\delta  , \rho/2) \mid \tu \delta \text{ is an irreducible genuine pseudospherical representations of  } \tu M\}$$ 
to be the set of irreducible  quotients of genuine pseudospherical principal series representations with infinitesimal character $\rho/2$. 
According to Table 1 in \cite{ABPTV}, the size of $\SH _{\rho/2}$ is either 1, 2 or 4, and we enumerate representations in $\SH_{\rho/2}$ as $\{Sh_i\}$. 
We have the following important properties for $\SH _{\rho/2}$. 

\begin{proposition} \label{Shi} Suppose we have the same setting as in Section \ref{sectionsmall}, and suppose that $G$ is simply laced and split. Then,  \\
(1) $\SH _{\rho/2} \subseteq \prod _{\rho/2} ^s (\widetilde{G})=  \prod _{\rho/2} ^{\Or} (\widetilde{G})$, where $\Or$ is defined in Section \ref{sectionsmall}.\\
(2) There is a bijection between $\SH _{\rho/2}$ and $\prod _g (Z(\widetilde{G}))$. Furthermore, there is a bijection between  this and $P/(2P+R)$, where $P$ and $R$ are the weight lattice and root lattice, respectively. Therefore, the central character of a  representation  from $\SH_{\rho/2}$ can be specified by an element in $P/(2P+R)$. \\
(3) There is a bijection between $\prod _g (Z(\widetilde{G}))$ and $\{\mathcal{O}_i\}$, where $\{\mathcal{O}_i\}$ is the set of real forms of $\mathcal{O}$, the complex associated variety of  representations from $\SH _{\rho/2}$. 
\begin{proof}

The infinitesimal character of each representation from $\SH _{\rho/2}$ is $\rho/2$ due to construction. Notice that for the parameter of any $Sh_i\in \SH_{\rho/2}$, all simple integral roots are not simple. However, by Lemma 3.1 in \cite{IC1}  for a simple integral root $\alpha$, there exists a positive root system $\Psi_{\alpha}\subseteq \triangle (\mathfrak g , \mathfrak a)$ such that $\Psi _{\alpha}$ contains the simple integral roots and $\alpha$ is simple for $\Psi_{\alpha}$. It turns out that  such $\alpha$ is nonintegral in the new system and is a nonparity real root, and hence it is  in the $\tau$-invariant (cf. \cite{IC1}, Theorem 4.12), which shows part (1). 

 It is proved in \cite{ABPTV} that $Z(\tu{M})=Z(\tu{G})$ and there is a one-to-one correspondence between $\prod _g (Z(\tu{M}))$ and $[2P^{\vee} \cap R^{\vee}] / 2R^{\vee}$, which is isomorphic to $P/(2P+R)$ when $G$ is simply laced and hence part (2) follows. 

 It can be observed from Table 1 in \cite{ABPTV} and Table 2 that when $G$ is simply laced and split, $\#\{\Or _i\}=|\prod _g (Z(\widetilde{G}))|$, and hence part (3) follows. 
\end{proof}
\end{proposition}

\begin{remark}\label{Shiorb}
Because of Proposition \ref{singlereal} and \ref{Shi}, we can attach to each $Sh_i$ a pair $(\chi_i, \mathcal{O} _i)$, where $\chi _i$ is the central character of $Sh _i$, and $\overline{\mathcal{O}_i  }=$AV$(Sh_i)$. Then  for each $\tau \in $Out$(G)$, $Sh_i ^{\tau}$ (see the notation defined in Section \ref{actionaut}) is associated to the pair $(\chi_i ^{\tau}, \mathcal{O}_i ^{\tau})$, (i.e. as $\tau$ permutes the central characters, it also permutes the real associated varieties).
\end{remark}

\subsection{Translation functors across a nonintegral wall} \label{wallcross} We recall some basic tools: cross-actions, Cayley and inverse Cayley transforms before starting to construct new representations. Most of the material in this section can be found in \cite{IC1}, \cite{AT} and \cite{RT}. Fix an infinitesimal character $\lambda$. (In our case, $\lambda =\rho/2$). In order to compute characters for nonlinear groups, we need a family of infinitesimal characters containing $\lambda$, denoted $\mathcal{F}(\lambda)$, which is defined as follows.  We define 
\begin{center}
%$\triangle=\triangle(\gc,\hc)$, $\hc \subseteq \gc$ is a fixed Cartan subalgebra,\\
$\triangle ^{+}=\{\alpha \in \triangle | \langle\lambda ,\alpha^{\vee}\rangle >0\}$,\\
$P=\{\mu  \in \mathfrak{h}^{\ast} \thinspace | \thinspace \langle\mu, \alpha ^{\vee} \rangle \in \mathbb{Z}  \text{ for } \alpha\in \triangle \}$, the integral weight lattice, \\
$W_P(\lambda)= \{  w\in W\thinspace | \thinspace w\lambda -\lambda \in P\}$.
\end{center}
Let $\mathcal{F}(\lambda)$ be a family of representatives of $(W\cdot\lambda +P)/P$ containing $\lambda$, and hence it is clear that $\mathcal{F}(\lambda)$ is indexed by $W/W_P(\lambda)$: if $\gamma \in \mathcal{F}(\lambda)$, then $\gamma =y\lambda$ modulo $P$ for some $y\in W$  which is unique modulo $W_{P}(\lambda)$. So we can write \begin{center}$\mathcal{F}(\lambda)=\{\gamma _y = y\lambda \thinspace | \thinspace y\in W/W_P(\lambda)\}$.\end{center} In particular, $\lambda=\gamma _1$.
There is an obvious action of $W$ on $\mathcal{F}(\lambda)$: $w\ast \gamma _{y} :=  w^{-1} (\gamma_y +\mu (y, w))=\gamma _{yw},$ by picking some $\mu (y, w)\in P$. We  fix once and for all integral wights $\mu(y,w)\in P$ satisfying the above conditions and we want to use them to define the following. First let $\alpha$ be a nonintegral simple root in $\triangle^+$, $s_{\alpha}$ be the corresponding simple reflection. Then we define:\\ (a) the nonintegral wall-crossing functors  $\psi _{\alpha}$  and $\phi_{\alpha}$, where    $\psi _{\alpha}(X) : = \psi _{\gamma_y} ^{\gamma_{ys_{\alpha}}}(X)$, a functor realizes an equivalence of categories between $\mathcal{HC}(\mathfrak{g}, K)_{  \gamma_y }$ and  $\mathcal{HC}(\mathfrak{g}, K)_ {\gamma_{ys_{\alpha}}} $; its inverse is $\phi_{\alpha}$(see \cite{V5}, Proposition 7.3.3);\\ (b) the cross action of $W$: let $\gamma=(H', \Gamma, \overline{\gamma})$ be a $(\gamma_y)$-regular character, $w\in W$, then  the regular character $w\times \gamma = (H', w\times \Gamma, w\times \overline{\gamma})$ is defined by $w\times\overline{ \gamma}=\overline{\gamma}+ \mu (y,w)$  and $w\times\Gamma=\Gamma\otimes \mu (y,w)\otimes \partial \rho(w)$, where  $ \partial \rho(w)  := w\cdot (\rho _i -2\rho_{ic})  - (\rho _i -2\rho_{ic})$, $\rho_i $ (resp. $\rho_{ic}$) denotes the half-sum of positive imaginary (resp. compact imaginary) roots that make $\overline{\gamma}$ dominant. Note that $w\times\overline{ \gamma}$ defines the same infinitesimal character as $\gamma_{yw}$.

\begin{remark}\label{infcharequiv}
Let $\alpha_1,\cdots,\alpha _p$ be simple roots and $s_1,\cdots, s_p$ be the corresponding reflections. If $w=s_p\cdots s_1 \in W_P(\lambda)$, we can define $\mu (1,w)=w\lambda-\lambda$, which is equal to  $\mu(1,s_1)+\mu(s_1,s_2)+ \mu (s_1s_2, s_3)+\cdots +\mu (s_1s_2\cdots s_{p-1} , s_p)$. Thus, $w\times \overline{\gamma} =\overline{\gamma} $ if and only if $w\in W_P(\lambda),$ where $\overline{\gamma} \sim \lambda.$ 
\end{remark}

We also need some basic facts about Cayley and inverse Cayley transforms. The related concepts can be found in various references (e.g. \cite{RT}, \cite{V5}). Here we just introduce some notation and quote some important facts.

    Let $\gamma= (H, \Gamma, \overline{ \gamma})$ be a $\lambda$-regular character. Assume $\alpha$ is a nonintegral root, then we can define Cayley (or inverse Cayley) transform on $\gamma$ (see Section 5 of \cite{RT}) through $\alpha$ if $\alpha$  is noncompact imaginary (real, resp.) and this action is denoted by $c^{\alpha} (\gamma)=\gamma^{\alpha}$ (or $c_{\alpha}(\gamma)=\gamma_{\alpha}$, resp.) Note that after  Cayley (inverse Cayley, resp.) transform, we get a new $\lambda$-regular character, say, $\gamma ^{\alpha}=(H^{\alpha},\Gamma^{\alpha},\overline{ \gamma^{\alpha}})$ (or $\gamma_{\alpha} =(H_{\alpha}, \Gamma _{\alpha},\overline{ \gamma_{\alpha}})$, resp.), which has infinitesimal character $\lambda$ and $I(\gamma ^{\alpha})$ (or $I(\gamma _{\alpha})$, resp.) has the same central character as the original representation $I(\gamma)$. For convenience, we call both operators $c^{\alpha}$ and $c_{\alpha}$  Cayley transforms through the root $\alpha$.

   Now we are ready to state the result of Vogan describing translation functors across a nonintegral wall.

\begin{theorem}\label{wallcrossing}
(\cite{IC1}, Corollary 4.8 and Lemma 4.9) Let $\gamma$ be a genuine $\lambda$-regular character of $G$. Suppose $\alpha$ is a nonintegral simple root in $\triangle ^{+}(\overline{\gamma})$. Then, with the translation functor $\psi_{\alpha}$ defined by the weight $\mu_{\alpha}$ fixed above, we have:\\
$\psi_{\alpha}(  J (\gamma) )= J(     (\gamma +\mu _{\alpha}) ^{\alpha}    ) =J  (    (s_{\alpha}  \times \gamma) ^{\alpha}  ) $ if $\alpha$ is noncompact imaginary, \\ $\psi_{\alpha}( J (\gamma) )=J (     (\gamma +\mu _{\alpha}) _{\alpha}    ) =J (    (s_{\alpha}  \times \gamma) _{\alpha}  ) $ if $\alpha$ is real satisfying the parity condition, \\ $\psi_{\alpha}(J(\gamma) )=J (    \gamma +\mu _{\alpha}   ) =J  (s_{\alpha}  \times \gamma ) $ otherwise. 
\end{theorem}
\begin{remark} \label{r1}
It can be observed from Theorem \ref{wallcrossing} that $\psi _{\alpha} (J(\gamma))$ and $J(s_{\alpha} \times \gamma)$ have the same infinitesimal character and central character.
\end{remark}

\subsection{Construction related to Dynkin Diagrams} \label{piRD}
Now we are ready to construct representations in $\prod_{\rho/2} ^s (\tu{G})$.  As described in \cite{P}, to the Dynkin diagram $D$ of $\widetilde{G}$, we attach a finite abelian group denoted by $R_{D}$ as follows. Let  $\prod$ be the set of simple roots. Define 

$R_D=\{S\subseteq \prod \thinspace|  \thinspace S$ is strongly orthogonal, so that any $\beta \notin S$ is adjacent to an even number of elements in $S\}$.

In Table 3, we list out the elements in $R_D$ for simply laced groups using Dynkin diagrams. Note that the root is in the element in $R_D$ if and only if the corresponding node is filled.

\begin{table} \label{DD}
\caption{}
\begin{center}
\begin{longtable}{|c|c|p{8.2cm}|c|}\hline 
Type & $n$&$R_D$&$|R_D| $\\\hline
$A_{n-1}$&$2m$&  \begin{picture}(50,15)
	\put(0,0){\circle{4}}
	\put(2,0){\line(1,0){16}}
	\put(20,0){\circle{4}}
	\put(22,0){\line(1,0){16}}
	\put(40,0){\circle{4}}
          \put(42,0){\line(1,0){8}}	
          \put(-2,4){\tiny{1}}
          \put(18,4){\tiny{2}}
           \put(38,4){\tiny{3}}
            \put(65,4){\tiny{$n-3$}}
              \put(87,4){\tiny{$n-2$}}
                \put(110,4){\tiny{$n-1$}}
          \put(51,-3){$\cdots$}
          \put(65,0){\line(1,0){8}}
          \put(75,0){\circle{4}}
          \put(77,0){\line(1,0){16}}
          \put(95,0){\circle{4}}
          \put(97,0){\line(1,0){16}}
          \put(115,0){\circle{4} }
          \put(147,0){\scriptsize{$\{\phi\}$} }
          \end{picture} 
          
            \begin{picture}(50,10)
	\put(0,-2){\circle*{4}}
	\put(2,-2){\line(1,0){16}}
	\put(20,-2){\circle{4}}
	\put(22,-2){\line(1,0){16}}
	\put(40,-2){\circle*{4}}
          \put(42,-2){\line(1,0){8}}	
          \put(-2,2){\tiny{1}}
          \put(18,2){\tiny{2}}
           \put(38,2){\tiny{3}}
            \put(65,2){\tiny{$n-3$}}
              \put(87,2){\tiny{$n-2$}}
                \put(110,2){\tiny{$n-1$}}
          \put(51,-5){$\cdots$}
          \put(65,-2){\line(1,0){8}}
          \put(75,-2){\circle*{4}}
          \put(77,-2){\line(1,0){16}}
          \put(95,-2){\circle{4}}
          \put(97,-2){\line(1,0){16}}
          \put(115,-2){\circle*{4} }
          \put(133,-2){\scriptsize{$\{\alpha_1, \alpha_3,\cdots, \alpha_{n-1}\}$} }
          \end{picture}

&2\\\hline
$A_{n-1}$&$2m+1$&   \begin{picture}(50,15)
	% Circles
	\put(0,0){\circle{4}}
	\put(2,0){\line(1,0){16}}
	\put(20,0){\circle{4}}
	\put(22,0){\line(1,0){16}}
	\put(40,0){\circle{4}}
          \put(42,0){\line(1,0){8}}	
          \put(-2,4){\tiny{1}}
          \put(18,4){\tiny{2}}
           \put(38,4){\tiny{3}}
            \put(65,4){\tiny{$n-3$}}
              \put(87,4){\tiny{$n-2$}}
                \put(110,4){\tiny{$n-1$}}
          \put(51,-3){$\cdots$}
          \put(65,0){\line(1,0){8}}
          \put(75,0){\circle{4}}
          \put(77,0){\line(1,0){16}}
          \put(95,0){\circle{4}}
          \put(97,0){\line(1,0){16}}
          \put(115,0){\circle{4} }
          \put(147,0){\scriptsize{$\{\phi\}$} }
          \end{picture} 
          
&1\\\hline
$D_n$&$2m$&   
\begin{picture}(50,20)
	\put(0,0){\circle{4}}
	\put(2,0){\line(1,0){16}}
	\put(20,0){\circle{4}}
	\put(22,0){\line(1,0){16}}
	\put(40,0){\circle{4}}
          \put(42,0){\line(1,0){8}}	
          \put(-2,4){\tiny{1}}
          \put(18,4){\tiny{2}}
           \put(38,4){\tiny{3}}
            \put(65,4){\tiny{$n-3$}}
          \put(51,-3){$\cdots$}
          \put(65,0){\line(1,0){8}}
          \put(75,0){\circle{4}}
          \put(77,0){\line(1,0){16}}
          \put(95,0){\circle{4}}
         \put(97,1){\line(1,1){10} }
          \put(108.2,11.9){\circle{4}}
         \put(108.2, -11.9){\circle{4}}    
         \put(97,-1){\line(1,-1){10}}
         \put(104, 15){\tiny{$n-1$}}
           \put(108, -17){\tiny{$n$}}
           \put(140, 0){\scriptsize{$\{\phi\}$  }}
\end{picture}
&4\\ && \begin{picture}(100,32)
	\put(0,0){\circle{4}}
	\put(2,0){\line(1,0){16}}
	\put(20,0){\circle{4}}
	\put(22,0){\line(1,0){16}}
	\put(40,0){\circle{4}}
          \put(42,0){\line(1,0){8}}	
          \put(-2,4){\tiny{1}}
          \put(18,4){\tiny{2}}
           \put(38,4){\tiny{3}}
            \put(65,4){\tiny{$n-3$}}
          \put(51,-3){$\cdots$}
          \put(65,0){\line(1,0){8}}
          \put(75,0){\circle{4}}
          \put(77,0){\line(1,0){16}}
          \put(95,0){\circle{4}}
         \put(97,1){\line(1,1){10} }
          \put(108.2,11.9){\circle*{4}}
         \put(108.2, -11.9){\circle*{4}}    
         \put(97,-1){\line(1,-1){10}}
         \put(104, 15){\tiny{$n-1$}}
           \put(108, -17){\tiny{$n$}}
           \put(140, 0){\scriptsize{$\{\alpha_{n-1}, \alpha_{n} \}$  }}
\end{picture}
& \\ && \begin{picture}(100,32)
	\put(0,0){\circle*{4}}
	\put(2,0){\line(1,0){16}}
	\put(20,0){\circle{4}}
	\put(22,0){\line(1,0){16}}
	\put(40,0){\circle*{4}}
          \put(42,0){\line(1,0){8}}	
          \put(-2,4){\tiny{1}}
          \put(18,4){\tiny{2}}
           \put(38,4){\tiny{3}}
            \put(65,4){\tiny{$n-3$}}
          \put(51,-3){$\cdots$}
          \put(65,0){\line(1,0){8}}
          \put(75,0){\circle*{4}}
          \put(77,0){\line(1,0){16}}
          \put(95,0){\circle{4}}
         \put(97,1){\line(1,1){10} }
          \put(108.2,11.9){\circle*{4}}
         \put(108.2, -11.9){\circle{4}}    
         \put(97,-1){\line(1,-1){10}}
         \put(104, 15){\tiny{$n-1$}}
           \put(108, -17){\tiny{$n$}}
           \put(140, 0){\scriptsize{$\{\alpha_1, \alpha_3,\cdots,\alpha_{n-3}, \alpha_{n-1}\}$  }}
\end{picture}
 & \\  &&\begin{picture}(100,32)
	\put(0,0){\circle*{4}}
	\put(2,0){\line(1,0){16}}
	\put(20,0){\circle{4}}
	\put(22,0){\line(1,0){16}}
	\put(40,0){\circle*{4}}
          \put(42,0){\line(1,0){8}}	
          \put(-2,4){\tiny{1}}
          \put(18,4){\tiny{2}}
           \put(38,4){\tiny{3}}
            \put(65,4){\tiny{$n-3$}}
          \put(51,-3){$\cdots$}
          \put(65,0){\line(1,0){8}}
          \put(75,0){\circle*{4}}
          \put(77,0){\line(1,0){16}}
          \put(95,0){\circle{4}}
         \put(97,1){\line(1,1){10} }
          \put(108.2,11.9){\circle{4}}
         \put(108.2, -11.9){\circle*{4}}    
         \put(97,-1){\line(1,-1){10}}
         \put(104, 15){\tiny{$n-1$}}
           \put(112, -15){\tiny{$n$}}
           \put(140, 0){\scriptsize{$  \{\alpha_1, \alpha_3,\cdots,\alpha_{n-3}, \alpha_{n}\} $  }}
\end{picture}
&\\ &&& \\
\hline
$D_n$&$2m+1$&  \begin{picture}(50,20)
	\put(0,0){\circle{4}}
	\put(2,0){\line(1,0){16}}
	\put(20,0){\circle{4}}
	\put(22,0){\line(1,0){16}}
	\put(40,0){\circle{4}}
          \put(42,0){\line(1,0){8}}	
          \put(-2,4){\tiny{1}}
          \put(18,4){\tiny{2}}
           \put(38,4){\tiny{3}}
            \put(65,4){\tiny{$n-3$}}
          \put(51,-3){$\cdots$}
          \put(65,0){\line(1,0){8}}
          \put(75,0){\circle{4}}
          \put(77,0){\line(1,0){16}}
          \put(95,0){\circle{4}}
         \put(97,1){\line(1,1){10} }
          \put(108.2,11.9){\circle{4}}
         \put(108.2, -11.9){\circle{4}}    
         \put(97,-1){\line(1,-1){10}}
         \put(104, 15){\tiny{$n-1$}}
           \put(108, -17){\tiny{$n$}}
           \put(140, 0){\scriptsize{$\{\phi\}$  }}
\end{picture}
&2\\ &&\begin{picture}(100,32)
	\put(0,0){\circle{4}}
	\put(2,0){\line(1,0){16}}
	\put(20,0){\circle{4}}
	\put(22,0){\line(1,0){16}}
	\put(40,0){\circle{4}}
          \put(42,0){\line(1,0){8}}	
          \put(-2,4){\tiny{1}}
          \put(18,4){\tiny{2}}
           \put(38,4){\tiny{3}}
            \put(65,4){\tiny{$n-3$}}
          \put(51,-3){$\cdots$}
          \put(65,0){\line(1,0){8}}
          \put(75,0){\circle{4}}
          \put(77,0){\line(1,0){16}}
          \put(95,0){\circle{4}}
         \put(97,1){\line(1,1){10} }
          \put(108.2,11.9){\circle*{4}}
         \put(108.2, -11.9){\circle*{4}}    
         \put(97,-1){\line(1,-1){10}}
         \put(104, 15){\tiny{$n-1$}}
           \put(111, -15){\tiny{$n$}}
           \put(140, 0){\scriptsize{$\{\alpha_{n-1},\alpha_n\}$  }}
\end{picture}
& \\ &&&
\\\hline
$E_6$ & &  \begin{picture}(50,31)
	\put(0,20){\circle{4}}
	\put(2,20){\line(1,0){16}}
	\put(20,20){\circle{4}}
	\put(22,20){\line(1,0){16}}
	\put(40,20){\circle{4}}
          \put(42,20){\line(1,0){16}}	
          \put(60,20){\circle{4}}
          \put(62,20){\line(1,0){16}}
          \put(80,20){\circle{4}}
          \put(-2,24){\tiny{1}}
          \put(18,24){\tiny{2}}
           \put(38,24){\tiny{3}}
           \put(58,24){\tiny{4}}
           \put(78,24){\tiny{5}}
           \put(40,18){\line(0,-1){15}}
           \put(40,2){\circle{4}}
           \put(42,-3){\tiny{6}}
            \put(140, 20){\scriptsize{$\{\phi\}$  }}
\end{picture}

&1
    \\\hline
$E_7$&&   \begin{picture}(50,10)
	\put(0,0){\circle{4}}
	\put(2,0){\line(1,0){16}}
	\put(20,0){\circle{4}}
	\put(22,0){\line(1,0){16}}
	\put(40,0){\circle{4}}
          \put(42,0){\line(1,0){16}}	
          \put(60,0){\circle{4}}
          \put(62,0){\line(1,0){16}}
          \put(80,0){\circle{4}}
          \put(82,0){\line(1,0){16}}
          \put(100,0){\circle{4}}
          \put(-2,4){\tiny{1}}
          \put(18,4){\tiny{2}}
           \put(38,4){\tiny{3}}
           \put(58,4){\tiny{4}}
           \put(78,4){\tiny{5}}
           \put(98,4){\tiny{6}}
           \put(60,-2){\line(0,-1){15}}
           \put(60,-18){\circle{4}}     
           \put(62,-23){\tiny{7}}
            \put(140, 0){\scriptsize{$\{\phi\}$  }}
\end{picture}
&2\\
&&
 \begin{picture}(50,40)
	\put(0,8){\circle*{4}}
	\put(2,8){\line(1,0){16}}
	\put(20,8){\circle{4}}
	\put(22,8){\line(1,0){16}}
	\put(40,8){\circle*{4}}
          \put(42,8){\line(1,0){16}}	
          \put(60,8){\circle{4}}
          \put(62,8){\line(1,0){16}}
          \put(80,8){\circle{4}}
          \put(82,8){\line(1,0){16}}
          \put(100,8){\circle{4}}
          \put(-2,12){\tiny{1}}
          \put(18,12){\tiny{2}}
           \put(38,12){\tiny{3}}
           \put(58,12){\tiny{4}}
           \put(78,12){\tiny{5}}
           \put(98,12){\tiny{6}}
           \put(60,6){\line(0,-1){15}}
           \put(60,-10){\circle*{4}}     
           \put(62,-15){\tiny{7}}
            \put(140, 8){\scriptsize{$\{\alpha_1, \alpha_3, \alpha_7\}$  }}
\end{picture}

&\\
\hline
$E_8$&&
  \begin{picture}(50,20)
	\put(0,5){\circle{4}}
	\put(2,5){\line(1,0){16}}
	\put(20,5){\circle{4}}
	\put(22,5){\line(1,0){16}}
	\put(40,5){\circle{4}}
          \put(42,5){\line(1,0){16}}	
          \put(60,5){\circle{4}}
          \put(62,5){\line(1,0){16}}
          \put(80,5){\circle{4}}
          \put(82,5){\line(1,0){16}}
          \put(100,5){\circle{4}}
          \put(102,5){\line(1,0){16}}
          \put(120,5){\circle{4}}
          \put(-2,9){\tiny{1}}
          \put(18,9){\tiny{2}}
           \put(38,9){\tiny{3}}
           \put(58,9){\tiny{4}}
           \put(78,9){\tiny{5}}
           \put(98,9){\tiny{6}}
           \put(118,9){\tiny{7}}
           \put(80,3){\line(0,-1){15}}
           \put(80,-13){\circle{4}}     
           \put(83,-15){\tiny{8}}
            \put(140, 5){\scriptsize{$\{\phi\}$  }}
\end{picture}

&1\\\hline
\end{longtable}
\end{center}
\end{table}

\begin{lemma} \label{PR}
There is a one-to-one correspondence between $R_D$ and $(Z(G_{\mathbb{C}})_2)^{\wedge}$, the characters of elements in $Z(G_{\mathbb{C}})$ of order 2. The latter  is isomorphic to $P/(2P+R)$ as group, and hence $R_D$ can be parametrized by the elements in $P/(2P+R)$. 

\begin{proof}
Denote $Z=Z(G_{\mathbb{C}})$ and $Z_2 =Z(G_{\mathbb{C}})_2$.

From the exact sequence 
\begin{center}
$1\to Z_2 \to  Z \to Z/Z_2\to 1 $,
\end{center}
we have another exact sequence \begin{center}
$1\to (Z/ Z_2)^{\wedge} \to  Z^{\wedge} \to Z_2^{\wedge}\to 1 $.
\end{center}
Notice that $Z^{\wedge}\simeq P/R$. Write $Z_2 =\{ exp(2\pi i \tau ^{\vee}) | \tau^{\vee} \in X_{\ast}\otimes \mathbb{C},  exp(2\pi i (2\tau ^{\vee})) =1 \}$.
Then $(Z/Z_2)^{\wedge}\simeq (2P+R)/R$, since $\gamma \in P$ such that $\gamma | _{Z_2}=1$ (i.e. $\gamma ( exp(2\pi i \tau ^{\vee})  )=exp (2\pi i \langle\gamma , \tau^{\vee}\rangle)  =1   $ ) if and only if $\langle\gamma , \tau ^{\vee}\rangle=1$, if and only if $\gamma \in 2P+R$.  Therefore, $Z_2 ^{\wedge} \simeq P/(2P+R)$ from the above exact sequence. 

 Associate  to each $S=\{\alpha_1, \cdots,  \alpha_p \}\in R_D$ an element $w_S=s_{\alpha_1}\cdots s_{\alpha_p} \in W$, then we have a map sending elements in $R_D$ to $P/(2P+R)$ by $S\to w_S (\rho/2)-\rho/2$. This is a bijection by counting the elements in $R_D$ and $P/(2P+R)$ case by case.
\end{proof}
\end{lemma}

\begin{remark} \label{pr}
When $G$ is simply laced,  $P/(2P+R) \cong  (2P^{\vee}\cap R^{\vee} ) /2R^{\vee}$, which parametrizes genuine central characters $\Pi _g (Z(\tu G))$ (see Section 3 in \cite{ABPTV}).

\end{remark}
We will show that we can get a subset of  representations in $\prod _{\rho/2} ^{s} (\widetilde{G})$ from each $Sh_i$ by a sequence of Cayley transforms or wall-crossings through the simple roots in $S\in R_D$.
 
 Associate  to each $S=\{\alpha_1, \cdots,  \alpha_p \}\in R_D$ an element $w_S=s_{\alpha_1}\cdots s_{\alpha_p} \in W$, and let  $c_S=c_{\alpha_1}\cdots c_{\alpha_p}$ and $\psi_S=\psi_{\alpha_1} \cdots \psi_{\alpha_p}$ be the corresponding Cayley transform and wall-crossing functor respectively. 
\begin{lemma}\label{wp}
For every $w_S$, $S\in R_D$, $w_S\in W_P(\rho/2)$, and hence $w_S\times \overline{\gamma}=\overline{\gamma}$, where $\gamma\sim\rho/2$, by  Remark \ref{infcharequiv}.
\begin{proof} Let $S=\{\alpha_1, \cdots,  \alpha_p \}\in R_D$. Then $w_S(\rho/2)-\rho/2=s_{\alpha_1}\cdots s_{\alpha_p}(\rho/2)-\rho/2=-\langle\rho/2,\alpha_1^{\vee}\ra\alpha_1-\cdots -\la\rho/2,\alpha_p^{\vee}\ra\alpha_p$.\\ For each simple root $\beta\notin S$, $\beta $ is adjacent to even numbers of $\alpha_i$'s, and hence $\la w_S(\rho/2)-\rho/2, \beta^{\vee}\ra\in \mathbb{Z}$. For $\beta=\alpha_i$ some $i$, $\la w_S(\rho/2)-\rho/2,\beta^{\vee}\ra=-\la\rho/2,\alpha_i^{\vee}\ra  \la\alpha_i,\alpha_i^{\vee}\ra\in \mathbb{Z}$. Therefore,  $w_S\in W_P(\rho/2)$. 
\end{proof}
\end{lemma}

\begin{lemma}\label{keyRD}
(1) $\psi _S (Sh _i)$ is in $ \prod _{\rho/2} ^s (\widetilde{G})$ for every $S\in R_D$ and $Sh_i\in \SH _{\rho/2}$;\\
(2) $\bigcup\limits _{Sh_i \in \SH _{\rho/2}}\{\psi _S(Sh_i) \thinspace | \thinspace S\in R_D  \}  = \bigcup\limits _{Sh_i \in \SH _{\rho/2}}\{ c _S(Sh_i) \thinspace | \thinspace S\in R_D  \} $.
%\begin{center}$\bigcup\limits _{Sh_i \in \SH}\{ c _S(Sh_i) \thinspace | \thinspace S\in R_D  \} \subseteq \prod _{\rho/2} ^s (\widetilde{G})$\end{center}.

\begin{proof}
Let $S\in R_D$ and $Sh_i\in \SH_{\rho/2}$. We apply $\psi_{S}$, a series of nonintegral wall-crossings $\psi_{\alpha}$'s, to $Sh_i$, and in each step, $\psi _{\alpha}(Sh_i)= P_{\gamma _y} ^{\gamma_{ys_{\alpha}}}(Sh_i\otimes F_{\mu (y,s_{\alpha})} )$, the projection of $Sh_i\otimes F_{\mu (y,s_{\alpha})} $ on to the Harish-Chandra modules with infinitesimal character $\gamma _{ys_{\alpha}  }$, where $\gamma _y$, $\gamma _{ys_{\alpha}}$, and $\mu (y,s_{\alpha})$ are described in  Section \ref{wallcross} . Note that $Sh_i\in \prod _{\rho/2} ^s  (\widetilde{G})=\prod _{\rho/2} ^{\mathcal{O}}  (\widetilde{G})$, we have $AV(I_{Sh_i})=\overline{\mathcal{O}}$ and hence $AV(I_{Sh_i\otimes F} ) = \overline{\mathcal{O}}$ for any finite dimensional $F$ by Lemma \ref{AV}. Therefore, $AV(I_{\psi_{\alpha} (Sh_i) })=\overline{\mathcal{O} }$ and  also $AV(I_{\psi_{S} (Sh_i)} )=\overline{\mathcal{O} }$. By Remark \ref{r1} and Lemma \ref{wp}, $\psi_{S} (Sh_i)$ has infinitesimal character $\rho/2$ and we conclude that $\psi_{S} (Sh_i)\in \prod_{\rho/2} ^{\Or} (\tu{G})=  \prod _{\rho/2} ^s (\widetilde{G})$.

For the second part of the proof, first we observe that all representations  on both sides have infinitesimal character $\rho/2$ due to (1) and the fact that the infinitesimal character doesn't change by Cayley transforms. Fix $S\in R_D$ and suppose that $S\neq \phi$. According to  Theorem \ref{wallcrossing}, each step of the wall-crossings in $\psi_S$  goes through the same  Cayley transform as in $c_S$, and hence $\psi _S (Sh_i)$  and $c_S (Sh_i)$ are specified by the same Cartan subgroup for every $Sh_i\in \SH _{\rho/2}$. It remains to show that there is a $Sh_j \in\SH _{\rho/2}$, $j\neq i$, such that $\psi _S( Sh_i)$ and $c_S (Sh_j)$ have the same central character.  Let $\gamma$ be the regular character of $Sh_i$.  The central character of $J(w_S \times \gamma)$ (as well as that of $\psi _S (Sh_i)$) can be specified by a nontrivial element $w_S(\rho/2)-\rho/2 \in P/(2P+R)$ by Lemma \ref{PR}. On the other hand, $P/(2P+R)$ also parametrizes $\prod _g (Z(\tu{G}))$ by Remark \ref{pr} and also the set $\{ c_S(Sh_k) \thinspace | \thinspace Sh_k \in \SH _{\rho/2} \}$, which means that there exists $j\neq i$ such that the central character of $c_S (Sh_j)$ is parametrized by $w_S(\rho/2)-\rho/2$ and hence $c_S (Sh_j) = \psi _S (Sh_i)$.
\end{proof}
\end{lemma}

 Denote the set we construct by 
  \begin{center}  $\prod_{R_D} (\widetilde{G}) =   \bigcup\limits _{Sh_i \in \SH}\{ c _S(Sh_i) \thinspace | \thinspace S\in R_D  \}$,
 \end{center}
then the following Theorem follows immediately from Lemma \ref{keyRD}.

\begin{theorem}
$\prod_{R_D} (\widetilde{G}) \subseteq  \prod _{\rho/2} ^s (\widetilde{G})$.
\end{theorem} 

\begin{remark} \label{number} If $|  \prod_g (Z(\tu{G}))  |=p,$ then $|\prod_{R_D} (\widetilde{G}) |=p^2$. Indeed, by Proposition \ref{Shi} (2), $|  \prod_g (Z(\tu{G})) | = | \mathcal{SH} _{\rho/2}|$. It can be observed from Table \ref{DD} that $|R_D|=|\mathcal{SH }_{\rho/2}|$ for each type, and hence $  |\prod_{R_D} (\widetilde{G}) | = |\mathcal{SH} _{\rho/2}| |R_D|=p^2$.
\end{remark}

\subsection{Exhaustion of $\prod _{\rho/2} ^s (\tu{G})$}
We have shown that $\prod _{R_D}(\widetilde{G})  \subseteq \prod _{\rho/2} ^s (\widetilde{G})$. In this section we will show by counting the elements in $ \prod _{\rho/2} ^s (\widetilde{G})$ that this is in fact an equality  when $G$ has type $A_{n-1}$ or $D_n$.

Fix a central character $\widetilde{\chi}$ of $\widetilde{G}$. Let $\prod _{\rho/2} ^s (\widetilde{G}) _{\widetilde{\chi}}$ be the subset of representations in $\prod _{\rho/2} ^s (\widetilde{G})$  with central character $\widetilde{\chi}$. The goal is to count $|  \prod _{\rho/2} ^s (\widetilde{G}) _{\widetilde{\chi}}  |$. Take a 
block $\mathcal{B}$ of representations with central  character $\widetilde{\chi}$ and infinitesimal character $\rho/2$, and consider $\prod(\rho/2), \triangle(\rho/2), W(\rho/2)$, the simple roots of the integral root system, the integral root system for $\rho/2$, and the integral Weyl group, respectively. Let $\mathbb{Z}[\mathcal{B}]$ be the $\mathbb{Z}$-span of the set of standard modules $I(\gamma _j)$, where each $\gamma _j$ is a $\rho/2$-regular character in $\mathcal{B}$. Then $W(\rho/2)$ acts on $\mathbb{Z}[\mathcal{B}] $ by the \textit{coherent continuation action} (see \cite{IC4} or \cite{parkcity}) and this action is denoted by  $w\cdot I(\gamma)$, or simply $w\cdot \gamma$ for $w\in W(\rho/2)$ and $\gamma \in \mathcal{B}$.

Consider $\{ J(\gamma) |\gamma \in \mathcal{B} \}$, the set of irreducible quotients of $\{ I(\gamma) |\gamma \in \mathcal{B} \}$, as another basis of $\mathbb{Z}[\mathcal{B}]$, we have

\begin{lemma} \label{taucoh}
Let $\alpha \in \prod (\rho/2)$, $\gamma \in \mathcal{B}$, then $s_{\alpha}\cdot J(\gamma)  = -J(\gamma)$ if and only if $\alpha \in \tau (J(\gamma))$.
\begin{proof}
Let $\alpha\in \prod (\rho/2)$ and let $\lambda$ be an infinitesimal character which is singular for $\alpha$. Define a coherent family with  $\pi (\rho/2) = J(\gamma)$. Then we have the identity \begin{center}$\pi (\rho/2) + \pi (s_{\alpha} (\rho/2 )) =\psi _{\lambda} ^{\rho/2} \circ \psi _{\rho/2} ^{\lambda} (\pi (\rho/2) ).$ \end{center}
Notice that $\alpha\in \tau(J(\gamma))$ if and only if $\psi _{\rho/2} ^{\lambda} (\pi (\rho/2) )=0$, which is equivalent to $\psi _{\lambda} ^{\rho/2} \circ \psi _{\rho/2} ^{\lambda} (\pi (\rho/2) )=0$, since the functor of  push-to or push-off walls  is injective.  We conclude that $\alpha\in \tau (J(\gamma))$ if and only if $J(\gamma) = \pi (\rho/2) =-  \pi (s_{\alpha} \rho/2) = -s_{\alpha} \cdot (J(\gamma))$ (the last equality holds by the definition of coherent continuation action). 
\end{proof}
\end{lemma}

\begin{lemma} \label{weylcalkey}
 $|  \prod _{\rho/2} ^s (\widetilde{G}) _{\widetilde{\chi}}  | =\dim Hom_{W(\rho/2)} (sgn, \mathbb{Z}[\mathcal{B}])$ 
\begin{proof}

Let $\pi = J(\gamma), \gamma \in \mathcal{B}$. Then by  Lemma \ref{taucoh}, $\pi \in  \prod _{\rho/2} ^s (\widetilde{G}) _{\widetilde{\chi}}$ if and only if $s_{\alpha} \cdot \pi =-\pi $ for all $\alpha \in \prod (\rho/2)$, which is equivalent to saying that $W(\rho/2)$  acts on $\pi$ as the sign representation. Thus the lemma follows.  
\end{proof}
\end{lemma}

Therefore, to count the left hand side, we just need to count the right hand side in  Lemma \ref{weylcalkey}. More precisely, we need to analyze the $W(\rho/2)$-representation $\mathbb{Z}[\mathcal{B}]$ in order to count the right hand side.

  The first observation is that  it makes counting more convenient if we consider a special block  $\mathcal{D}$, which is equivalent to  $\mathcal{B}$, instead of $\mathcal{B}$, but with an infinitesimal character other than $\rho/2$, and in a different chamber.
  
  The following lemma tells what this block is.
\begin{lemma} \label{simple}
Let $\rho = \rho(\triangle^+ )$, and $\lambda _0 = \rho/2$. Let $\prod$ be the simple roots in the chamber of $\rho/2$. Then we can find $w\in W$ with the following properties. Let $(\triangle ')^+ = w \triangle ^+$ and $\prod ' = w\prod  =\{\alpha _i\}$. There is some $\alpha_k \in \prod '$ such that if setting $\lambda = w\rho -\frac{1}{2} \lambda _k ^{\vee} ((\triangle ')^+)$, where $\lambda _ k ^{\vee}$ is the fundamental weight for $\alpha _k$, we have (1) $\la\lambda , \alpha _ k ^{\vee}\ra= 1/2 $ and $\la\lambda , \alpha _j ^{\vee}\ra = 1$ elsewhere; (2) $\triangle (\lambda) =\triangle (\lambda _0)$, (and $\prod (\lambda) = \prod (\lambda _0)$ as well).  Therefore, for type $A_{n-1},$ we can always move to a block $\mathcal{D}$ (through nonintegral wall-crossing equivalence) with infinitesimal character $\lambda_{\mathcal{D} } =\lambda$, such that every root in $\prod (\lambda)$ is simple for the entire root system; for type $D_n, n \geq 4,  E_6, E_7, E_8$, we can move to a block $\mathcal{D}$ with infinitesimal character $\lambda$, such that every root in $\prod (\lambda)$ is simple for the entire root system but one, say, $\alpha$.

\begin{proof}
This lemma is proved by a case by case calculation.
\end{proof}
\end{lemma}

For example, consider type $D_4$, $G=Spin(4,4)$. In this case, $\rho/2 = (\frac{3}{2}, 1, \frac{1}{2},0 )$, $\prod  = \{ e_1-e_2, e_2-e_3, e_3\pm e_4 \}$, and $\prod (\rho/2) = \{ e_1\pm e_3, e_2\pm e_4\}$. The infinitesimal character $\lambda$ to be chosen in Lemma \ref{simple} is $\lambda =(\frac{5}{2}, 1, \frac{3}{2}, 0 )$ and $\prod '$, the simple roots in this chamber, is $\{  e_1-e_3, e_3-e_2, e_2\pm e_4\}$, and it can be easily seen that $\prod (\lambda) =\prod (\rho/2) $ and the only simple integral root that is not simple in the entire root system is $\alpha= e_1+e_3$.

Due to the equivalence of block $\mathcal{B}$ and block $\mathcal{D}$, we'll focus on analyzing $\mathbb{Z}[\mathcal{D}]$ from now on and then count dim$_{W(\lambda)} (sgn, Z[\mathcal{D}])$. (Notice that from now on, $\lambda$ stands for the infinitesimal $\lambda _{\mathcal D}$ specified in Lemma \ref{simple}.)

Now take a closer look at the coherent continuation action of $W(\lambda)$ on $\mathbb{Z}[\mathcal{D}]$. The formulas of the coherent continuation action can be derived from the formula of the Hecke operators (see \cite{RT}, Definition 9.4). More precisely, let $\beta$ be a simple root, consider $T_{s_{\beta}}$, which we denote $T_{\beta}$ for simplicity, an operator acting on $\Z [q,q^{-1}] [\mathcal{D}]$ (the formulas are  given in \cite{RT}), and the coherent continuation action of $s_{\beta} \in W$ can be defined as 
\begin{equation}\label{hecke}
s_{\beta} \cdot  \gamma = - T_{\beta} (1) (\gamma) \quad\text{ with each term $\delta$ multiplied by $(-1) ^{l (\gamma) -l (\delta)}$,  }
\end{equation}
where  $\gamma \in \mathcal{D}$ and $l$ is a length function defined on parameters and it can be looked up in \cite{parkcity}. Notice that in (\ref{hecke}), $\beta$ could be nonintegral and  hence $s_{\beta} \cdot  \gamma$ is possibly outside of $\Z [\mathcal{D}]$. But if $\beta$ is integral, $s_{\beta}\in W(\lambda)$ and the following proposition gives explicit formulas for the action of $W(\lambda)$ on $\Z [\mathcal{D}]$, which can be derived from Definition 9.4 of \cite{RT} and (\ref{hecke}).

\begin{proposition} \label{cohcon}
Fix $\gamma \in \mathcal{D}$ and $\beta\in \prod (\lambda)$. Furthermore, suppose $\beta$ is simple in the whole root system.  Let $s:= s_{\beta}\in W(\lambda)$.\\
(a) If $\beta$ is  complex or real for $\gamma$, then $s\cdot \gamma=s\times \gamma.$\\
(b) If $\beta$ is compact imaginary for $\gamma$, then $s \cdot \gamma=-\gamma$.\\
(c) If $\beta$ is noncompact imaginary for $\gamma$, then $s\cdot \gamma =- s\times \gamma+c _{\beta}(\gamma)$.
\end{proposition}

The formulas in Proposition \ref{cohcon} are the same as those for the linear case (or say, when the infinitesimal character is integral) (see Theorem 4.12 in \cite{IC1}). Notice that in Proposition \ref{cohcon}, $\beta$ is simple for the whole root system; if $\beta =\alpha$, the only simple integral root which is not simple for the whole root system (see Lemma \ref{simple}), we will need an additional formula for $s_{\alpha}\cdot \gamma$ and this will be discussed later.

From Proposition \ref{cohcon}, we can see that the coherent continuation action is closely related to the cross action, so we also consider the cross action of $W(\lambda)$ on $\mathcal{D}$. Notice that two $\lambda$-regular characters $\gamma _i =(\widetilde{H_i}, \Gamma _i, \overline{\gamma_i})$ and $\gamma _j =(\widetilde{H_j}, \Gamma _j, \overline{\gamma_j})$ from $\mathcal{D}$ are in the same cross action orbit if and only if $\widetilde{H_i} = \widetilde{H_j}$.  Indeed, if $\widetilde{H_i} = \widetilde{H_j} = \widetilde{H}$, then $\Gamma _i$ and $\Gamma _j$ agree on $Z(\widetilde{G})$, since $Z (\widetilde{H}) =Z(\widetilde{G}) \widetilde{H_0}$ (by Proposition \ref{center} (1)). Since $\gamma_i$ and $\gamma _j$ are in the same block, $\overline{\gamma_i}$ and $\overline{\gamma _j}$ define the same infinitesimal  character, say, $\overline{\gamma _i}\sim \overline{\gamma _j} \sim \lambda$ and hence $\gamma _j = w\times \gamma _i$ for some $w\in W(\lambda)$.  
Enumerate the Cartan subgroups of $\widetilde{G}$ as $\{ \widetilde{H_1}, \cdots , \widetilde{H_l}\}$, and pick a $\lambda$-regular character  $\gamma _j$ specified by $\widetilde{H_j}$, then $\{\gamma _1, \cdots, \gamma _l \}$ is a set of representatives of the cross action orbits of $W(\lambda)$ on  $\mathbb{Z}[\mathcal{D}]$. 

Let $W_{\gamma _j} = \{ w \in W(\lambda) \thinspace | \thinspace w\times \gamma_j =\gamma _j  \}$ be the cross stabilizer of $\gamma _j $ in $W(\lambda)$. Then we have the following proposition.

\begin{proposition}  \label{decomp}
$\mathbb{Z}[\mathcal{D}] \simeq \oplus _j  Ind ^{W(\lambda)} _{W_{\gamma _j}   }  (\epsilon _j)$, where $\epsilon _j$ is a one-dimensional representation of $W_{\gamma _j}$ such that for $w\in W_{\gamma_j}$, $w \cdot \gamma _j = \epsilon _j (w) \gamma _j +$ other terms from more split Cartan subgroups.
\begin{proof}
The proof is similar to the linear case (cf. \cite{BV2}) by using formulas in Definition 9.4 of \cite{RT} and (\ref{hecke}).
\end{proof}
%\begin{proof}
%Since $W_{\gamma_j}$ is generated by $\{ s_{\beta}  | \beta \in \prod (\rho/2)\}$, it suffices to show that $s_{\beta}  \cdot I(\gamma _j) =\pm I(s_{\beta} \times \gamma _j )+$ other terms from more split Cartan subgroups. This is clear when $\beta$ is simple for the whole root system by Proposition \ref{cohcon}.

  %Consider $\alpha\in \prod(\rho/2)$, which is not simple (as listed in the table). Let $T_{\alpha}$ be the corresponding Hecke operator, the $\mathbb{Z}[q^{1/2}, q^{-1/2}]$-linear map from $\mathbb{Z}[\mathcal{D}] [q^{1/2}, q^{-1/2}]$ to itself (defined in \cite{parkcity}).  We can decompose $T_{\alpha} = T_{\alpha_1}T_{\alpha_2}\cdots T_{\alpha_m}$, where $\alpha_j$'s are simple. Here the $\alpha$'s are allowed to be non-integral. Notice that $$T_{\alpha}= -\phi _{\alpha} \psi _{\alpha} +q,$$ up to a sign,  where $\psi _{\alpha}$ and $\phi_{\alpha}$ are the functors of push-to and push-off walls, respectively, and also, $$\phi _{\alpha } \psi_{\alpha}  (I(\gamma)) =I(\gamma) +s_{\alpha} \cdot I(\gamma).$$ We conclude that $s_{\alpha} \cdot I(\gamma)=-T_{\alpha} (1) (\gamma)$, up to a sign. Using definition 9.4 in \cite{RT}, each $T_{\alpha_j } (I(\gamma  ))$ can be calculated explicitly, and so can $T_{\alpha} (I(\gamma))$. Therefore, it is not hard to see that $s_{\alpha}  \cdot I(\gamma _j) =\pm I(s_{\alpha} \times \gamma _j )+$ other terms from more split Cartan subgroups.
%\end{proof}
\end{proposition}

By Proposition \ref{decomp} and Frobenius reciprocity, the multiplicity of $sgn_{W(\lambda)}$ in $\mathbb{Z} (\mathcal{D})$ is  $[sgn_{W(\lambda)}:   \mathbb{Z}[\mathcal{D}]] =[sgn_{W(\lambda)} | _{W_{\gamma_j}} : \epsilon _j   ]$, which is equal to 0 or 1, since $sgn_{W(\lambda)} | _{W_{\gamma_j}} $ is one-dimensional. This means that we have reduced our goal to count the number of $\gamma _j$'s making $[sgn_{W(\lambda)} | _{W_{\gamma_j}} : \epsilon _j   ] = 1$, which is called condition $(\ast)$. Equivalently, $\gamma _j$ satisfies condition $(\ast)$  if 
\begin{center} 
\hspace*{5em}$sgn_{W(\lambda)}  | _{W _{\gamma _j}}  =\epsilon _j $     \hspace{10em} $(\ast)$
\end{center}

To reach this goal, we need to analyze $\epsilon _j$ for each $j$. By \cite{AT},  \begin{center}$W_{\gamma_j}=W^C (\overline{\gamma _j}) ^{\theta} \ltimes (W^i (\overline{\gamma _j}) \times W^r (\overline{\gamma _j}))$.\end{center}
So we can decompose $\epsilon _j = \epsilon ^C _j \otimes  \epsilon ^i _j \otimes  \epsilon ^r _j$, where $\epsilon ^C _j ,  \epsilon ^i _j , \epsilon ^r _j$ are characters of  $W^C (\overline{\gamma _j}) ^{\theta},$ $W^i (\overline{\gamma _j}),$ $W^r (\overline{\gamma _j})$, respectively. Notice that in the linear case (or say, when the block considered has integral infinitesimal character), we have $\epsilon _j = sgn _i$ for all $j$ (see \cite{BV2}), where $sgn _i = sgn _{  W^i (\overline{\gamma _j}) }$.  When $\tu{G}$ has type $A_{n-1}$, we are in a chamber where every $\beta\in \prod (\lambda)$ is simple, and then using formulas in Proposition \ref{cohcon}, an argument analogous to the linear case  shows that $\epsilon _j = sgn_i$ for every $j$. In fact, there is no compact imaginary root in the case of type $A_{n-1}$, and hence  the factor $W^i ( \Br{\gamma _j}) $ in $W_{\gamma_j}$ is trivial and therefore $sgn_i =1$. The following proposition follows.

\begin{proposition} \label{Aonedim}
 If $\widetilde{G}$ has type $A_n$, $\epsilon _j = 1 $ for all $\gamma _j$.
\end{proposition}
The following lemma is an easy result of Proposition \ref{Aonedim}.
\begin{lemma} \label{Aruleout}
For type $A_{n-1}$, $\gamma _j$ satisfies condition $(\ast)$  if and only if  there are no real integral roots for $\gamma _j$.

\end{lemma}

Analyzing $\epsilon _j$ for type $D_n$, $n\geq 4$, requires more work.  The formulas in Proposition \ref{cohcon} are not enough since there is an integral root $\alpha\in\prod (\lambda)$ which is not simple. So the first goal is to calculate $s_{\alpha}\cdot \gamma$, $\gamma\in\mathcal{D}$.  Notice that 
\begin{eqnarray*}
\prod (\lambda)=\prod (\rho/2) &=&\{ e_i -e_{i+2}, \thinspace {\scriptstyle 1 \leq i \leq n-2}, \thinspace e_{n-3}+e_{n-1}, e_{n-2}+e_n \}\\
\prod  &=& 
\begin{cases}
\{ e_i -e_{i+2}, \thinspace {\scriptstyle 1 \leq i \leq n-2}, \thinspace e_{n-1}-e_2, e_{n-2}+e_n \} & \text{ if $n$ is even }\\
\{ e_i -e_{i+2}, \thinspace {\scriptstyle 1 \leq i \leq n-2}, \thinspace e_n-e_2, e_{n-3}+e_{n-1} \} & \text{ if $n$ is odd } 
\end{cases}   \\
\text{Terefore, } \quad\alpha &=&
\begin{cases}
e_{n-3}+e_{n-1} &  \text{ if $n$ is even }\\
e_{n-2}+e_n &  \text{ if $n$ is odd } 
\end{cases}
\end{eqnarray*}

We decompose $s_{\alpha}= s_{\alpha_1} s_{\alpha_2} \cdots s_{\alpha_m}$, a product of simple reflections. 
\begin{equation}\label{simrefl}
s_{\alpha}=
\begin{cases}
 \nu s_{\scriptscriptstyle n-2,n}s_{\scriptscriptstyle\Br{\scriptscriptstyle n-2,n}}  s_{\scriptscriptstyle n-1,n-3} \nu  s_{\scriptscriptstyle n-1,n-3} s_{\scriptscriptstyle\Br{\scriptscriptstyle n-2,n}} s_{\scriptscriptstyle n-2,n} \nu &  \text{ if $n$ is even }\\
  s_{ \scriptscriptstyle n-2, n} \mu  s_{\scriptscriptstyle\scriptscriptstyle n-2,n}  s_{\scriptscriptstyle n-1,n-3} s_{\scriptscriptstyle \Br{\scriptscriptstyle n-1,n-3}} \mu  s_{\scriptscriptstyle \Br{\scriptscriptstyle n-3,n-1}} s_{\scriptscriptstyle n-1,n-3}s_{\scriptscriptstyle n-2,n} \mu  s_{\scriptscriptstyle n-2,n}&  \text{ if $n$ is odd}, 
\end{cases}
\end{equation}
where $s_{i,j}=s_{e_i-e_j}$ and $s_{\Br{i,j}}= s_{e_i+e_j}$, and $\nu = s_{n-4, n-2}\cdots s_{2,4} s_{2, n-1} s_{2,4} \cdots s_{n-4, n-2}$,
$\mu=  s_{n-3,n-1}s_{n-5, n-3}\cdots s_{2,4} \cdots s_{n-5,n-3} s_{n-3,n-1} s_{n-5, n-3} \cdots s_{2,4}\cdots s_{n-3,n-1}$.

Given $\gamma \in\mathcal{D}$. By (\ref{hecke}), to calculate $s_{\alpha}\cdot\gamma$, we have to calculate $T_{\alpha}(\gamma)$ first. Consider the decomposition $T_{\alpha} (\gamma) = T_{\alpha_1}T_{\alpha_2}\cdots T_{\alpha_m} (\gamma)$. Note that on the right hand side, the Hecke operation is calculated step by step. In each step, we have to deal with some $T_{\alpha_k} (\delta)$, where $\delta$ is the parameter of a standard module not necessarily belonging to block $\mathcal{D}$.  In fact, this is an "abstract" Hecke operation, and it should be denoted by $T_{\alpha_k}\cdot _{a} (\delta)$.  Taking an inner automorphism  $\phi _k$ of $\mathfrak{g}$ sending $(\lambda, \mathfrak{h}^{\ast})$ to  $(\overline{\delta}, \mathfrak{h}_{\delta} ^{\ast})$, we define   $T_{\alpha_k}\cdot _{a} (\delta) := T_{\phi _k  (\alpha_k)} (\delta)$. Here $\phi _k (\alpha_k)$ is a simple root in the chamber of $\delta$, and hence we can use the formulas in Definition 9.4 of \cite{RT} to calculate $T_{\phi _k(\alpha_k)} (\delta)$ in each step. 
\begin{eqnarray*} 
T_{\alpha} (\gamma) &=& T_{\alpha_1} \cdot _a T_{\alpha_2}\cdot _a \cdots \cdot _aT_{\alpha_m} (\gamma) \\
&= &       ( T_{\phi _1 (\alpha_1)}  (T_{\phi _2 (\alpha _2)}( \cdots    (  T_{\phi _m(\alpha_m)} (\gamma) ))\cdots ) \\
&=& p_1(q) \cdots p_m (q)  \phi _1 (\alpha_1)\times (\phi _2 (\alpha _2 \times \cdots (\phi _m (\alpha _m )\times \gamma ) ) \\ 
&&+ ( \text{terms from more split Cartan subgroups}),
\end{eqnarray*} 
where   $p_j (q) \in \mathbb{Z}[q, q^{-1}]$.
 According to (\ref{hecke}) and \cite{parkcity}, it turns out that at each step, we may define 
$$s_{\alpha_k} \cdot \delta = -T_{  \phi _k(\alpha_k )} (1) (\delta), \text{ if } \phi _k (\alpha_k) \text{ is real or imaginay for }  \delta$$ and 
$$s_{\alpha_k} \cdot \delta =T_{  \phi _k(\alpha_k )} (1) (\delta), \text{ if }  \phi _k (\alpha_k) \text{ is  complex for } \delta.$$
Let $t_{\gamma}$ be the number of occurrences imaginary roots in $\{\phi _j (\alpha_j) , 1\leqslant  j \leqslant m  \}$. An easy calculation shows that 
\begin{equation} \label{key}
s_{\alpha}\cdot \gamma = (-1)^{t_{\gamma}} s_{\alpha}\times \gamma +\text{(terms from more split Cartan subgroups)}.
\end{equation}

Notice that $\{\alpha_j\}_{j=1} ^m$ can be read off from (\ref{simrefl}). 

%When $n$ is even,  $m=9$ and \\
%$\{  \phi_j(\alpha_j) \} = \{   e_{n-3} +e_{n-2}, e_{n-3}+e_n, e_{n-3}-e_n, e_{n-2}+e_{n-1},e_{n-3}+e_{n-1},\\ e_{n-3}-e_{n-2}, e_{n-1}+e_n, e_{n-1}-e_n, e_{n-1}-e_{n-2} \}$.

%When $n$ is odd, $m=11$ and \\
%$\{  \phi_j(\alpha_j) \}  = \{  e_{n}-e_{n-2}, e_{n-3}+e_n, e_{n-3}+e_{n-2}, e_{n-1} +e_n, e_n-e_{n-1}, e_{n-2}+e_n, \\
%e_{n-2}+e_{n-1}, e_{n-2}-e_{n-1}, e_n-e_{n-3}, e_{n-2}-e_{n-3}, e_{n-2}-e_n \}$.

\begin{example}\label{funny}
Consider type $D_4$ and $\tu{G}=\tu{Spin}(4,4)$. In this case,  $\lambda =  (\frac{5}{2}, 1, \frac{3}{2}, 0  ),$ $\alpha = e_1+e_3$. According to Section \ref{wallcross}, we need a family of infinitesimal characters to define the cross action of $W$ on $\mca{D}$ and it  could be chosen to be $F(\lambda)=\{ (\frac{5}{2}, 1, \frac{3}{2}, 0  ), (\frac{3}{2}, 1,\frac{1}{2}, 0) , (\frac{3}{2}, 2, \frac{1}{2} , 0 )  \}$. Let $\gamma\in\mca{D}$ be parametrized by $( \underline{\frac{5}{2} \thinspace 1}\thinspace \frac{3}{2}, 0)$, meaning that $e_1-e_2$ is imaginary for $\gamma$ and $e_1+e_2, e_3\pm e_4$ are real for $\gamma$. (In other words, let $\gamma_0$ be the 
parameter for the pseudospherical principal series Ind$_{\tu M AN} ^{\tu G} (\tu \delta \otimes e^{ (\frac{5}{2}, 1, \frac{3}{2},0 )} \otimes 1)$, then $\gamma =c_{e_1-e_2}(\gamma_0)$ (see notation above Theorem \ref{wallcrossing}), which is obtained from $\gamma_0$ by taking Cayley transform through the root $e_1-e_2$). We apply a sequence of cross actions to $\gamma$ through the roots  $\{\phi _9(\alpha_9) ,  \cdots, \phi _1 (\alpha_1)\}$:

\begin{align*}
\gamma=&( \underline{\frac{5}{2} \thinspace 1}\thinspace \frac{3}{2}, 0)    \xrightarrow{ s_{e_2-e_3}  \times} ( \ul{ \frac{3}{2} \ths 1 } \ths \frac{1}{2} , 0)  
\xrightarrow{   s_{e_3-e_4 }\times }( \ul{2\ths\frac{3}{2} } \ths 0, \frac{1}{2}) \xrightarrow{ s_{e_3+e_4} \times} ( \ul{ \frac{3}{2} \ths 1 } \ths -\frac{1}{2} , 0) \xrightarrow{   s_{e_1-e_2 }\times } 
\\
&( \ul{ \frac{3}{2} \ths 2 } \ths -\frac{1}{2} , 0)  \xrightarrow{   s_{e_1+e_3 }\times } ( \ul{ \frac{1}{2} \ths 2 } \ths -\frac{3}{2} , 0) \xrightarrow{ s_{e_2+e_3} \times}( \ul{ \frac{1}{2} \ths 1 } \ths -\frac{3}{2} , 0) \xrightarrow{ s_{e_1-e_4} \times}  (\ul{0\ths\frac{3}{2} } \ths -2, \frac{1}{2})\\ &\xrightarrow{ s_{e_1+e_4} \times} ( \ul{ -\frac{1}{2} \ths 1 } \ths -\frac{3}{2} , 0) 
\xrightarrow{ s_{e_1+e_2} \times} ( \ul{ -\frac{3}{2} \ths 1 } \ths -\frac{5}{2} , 0) =s_{\alpha}\times\gamma
\end{align*}
Notice that we vary the chambers and fix the types of roots at each step. It can be seen that $\phi_6(\alpha_6) =e_1-e_2$ is the only imaginary root  for $\gamma$ among $\{  \phi_j(\alpha_j) \}$ and hence $t_{\gamma} =1$. We conclude that \begin{equation*} s_{\alpha}\cdot \gamma=-(s_{\alpha} \times \gamma) +\text{ (terms from more split Cartan subgroups)} \end{equation*} and hence $\epsilon (s_{\alpha} )=-1$.
 
\end{example}
Due to the remark above Proposition \ref{decomp}, we can choose each $\gamma_j$ properly and calculate the $\epsilon _j$'s according to the chosen $\gamma_j$'s. In fact, our goal is to rule out $\gamma_j$'s satisfying either of the following conditions.

\begin{itemize}
  \item (R) If there is a real integral root, then choose $\gamma_j$ such that $\alpha$ is real for $\gamma_j$.
  \item (C) If there are no real integral roots, and there is an orthogonal set of 4 nonintegral roots of the form $\{ e_p \pm e_q. e_r \pm e_s\}$, where $e_p \pm e_q$ are both imaginary, or both real, and one of $\{e_r\pm e_s\}$ is real, whereas the other is imaginary, then choose $\gamma _j$ such that  $\{  e_{n-3}\pm e_{n-2} , e_{n-1}\pm e_n\}$ is the desired quadruple. In this case $\alpha$ is a complex root. 
 
\end{itemize}
\begin{remark}
The $\gamma$ chosen in Example \ref{funny} is a parameter satisfying condition (C).
\end{remark}
With the setting, we have the following key lemma.

\begin{lemma} \label{Druleout}
Suppose that $\widetilde{G}$ has type $D_n$,  $n\geqslant 4$.  If $\gamma _j$ satisfies condition (R) or (C) then $\gamma _j$ doesn't satisfy condition $(\ast)$.

\begin{proof}
To show that the chosen $\gamma _j$ fails to satisfy condition $(\ast)$, we will pick a $w\in W_{\gamma _j}$ and show that $\epsilon _j(w)$ and $sgn (w)$ do not coincide.  In either case, we have to calculate $\epsilon _j (s_{\alpha})$ for $\gamma _j$. By (\ref{key}), we just need to count the number $t_{\gamma_j}$ for the chosen $\gamma_j$.

Suppose that $n$ is even. If  $\gamma _j$ satisfies condition (R), $e_{n-3}-e_{n-1}$ is also a real integral root, and the roots  $e_{n-2}\pm e_{n-1}, e_{n-2}\pm e_{n-3}, e_n\pm e_{n-1}, e_n\pm e_{n-3}$ can be arranged so that each of them is either real or complex. Therefore, $t_{\gamma_j}=0$, and hence $\epsilon _j (s_{\alpha})=1$. This result follows for the odd case  by applying the same argument.   Since $s_{\alpha} \in W_{\gamma _j}$ and sgn$(s_{\alpha})=-1$, $\gamma _j$ fails to satisfy condition $(\ast)$.

Now suppose that $\gamma _j$ does not satisfy condition (R), but satisfies condition (C). It can be counted that $t_{\gamma _j} =1$ or 3 (see Example \ref{funny} for the calculation), which implies $\epsilon _j (s_{\alpha}) = -1$. Let $w= s_{  e_{n-3 } -e_{n-1}  } s_{  e_{n-3 } +e_{n-1}  }   s_{  e_{n-2 } -e_{n}  } s_{  e_{n-2 } +e_{n}  }$. Notice that  $w\in W_{\gamma _j}$. %Indeed, the type of  roots $  e_{n-3}\pm e_{n-2} , e_{n-1}\pm e_n$ remain the same for $w\times \gamma _j$ and hence $w\in W_{\gamma_j}$.  
When $n$ is even (odd, respectively), we have $ e_{n-3 } -e_{n-1},  e_{n-2 }  \pm e_{n} $ ($ e_{n-3 } \pm e_{n-1},  e_{n-2 } -e_{n}$, respectively) are simple and complex, so $\epsilon _j (s_{\beta})  =1$ for every $\beta$ from these three root, and hence $\epsilon _j (w) = \epsilon_j (s_{\alpha} )=-1$. But it is obvious that $sgn(w)=1$. Therefore $\gamma _j$ fails to satisfy condition $(\ast)$.  
\end{proof}
\end{lemma}
\begin{theorem} \label{exhaustAD}
For type $A_{n-1}$ and $D_n$, $n\geqslant 4$, we have  $\prod _{\lambda} ^s (\widetilde{G})= \prod _{R_D}(\widetilde{G}) $.
\begin{proof}
It suffices to count the number of $\gamma_j$'s in Proposition \ref{decomp} satisfying condition $(\ast)$. Fixing a genuine central character, we shall show that this number is 2 (1, respectively) for type $A_{n-1}$, when $n$ is even (odd, respectively), and it is 4 (2, respectively) for type $D_n$, when $n$ is even (odd, respectively). We will give the proof for the even case only, and a similar argument applies for the odd case. 

For type $A_{n-1}$, we claim that if the real rank of the Cartan subgroup $H_j$ is at least $n/2$ (when $n$ is even) or $(n-1)/2$ (when $n$ is odd) , then there exists a real integral root for $\gamma _j$, and hence such $\gamma _j$ can be ruled out by Lemma \ref{Aruleout}. 

When $n$ is even, we enumerate all Cartan subgroups (on the level of linear groups) as $\{ H_{n/2 -1}, H_{n/2}, \cdots, H_{n-2}, H_{n-1} \}$, where  the real rank of $H_j$ is $j$. Let $\gamma _{n-1}$ be the parameter of the principal series, and $\alpha_k = e_{2k-1}-e_{2k}$, $1 \leqslant k \leqslant n/2$, then we pick $\gamma _  {n-1-k}= c_{\alpha _k}\cdots c_{\alpha _2} c_{\alpha_1} (\gamma _{n-1})$ to be the representative of the cross action orbit specified by $H_{n-1-k}$, $1 \leqslant k \leqslant n/2$. Notice that when $k\leqslant n/2 -2$, $e_{n-2}-e_n$ is a real integral root for $\gamma _{n-1-k}$, which means that we can rule out $\gamma _j$, for $n/2 +1\leqslant  j\leqslant n-1$. Only $\gamma _{n/2-1}$ and $\gamma _{n/2}$ are not ruled out, and they are exactly the $\gamma _j$'s  satisfying condition $(\ast)$ since the number of $\prod _{R_D} (\widetilde{G})$ with a fixed central character is also 2. Hence the theorem follows for type $A_{n-1}$, when $n$ is even.  When $n$ is odd, it can be shown by a similar argument that the only $\gamma_j$ satisfying condition $(\ast)$ comes from the fundamental Cartan.

For type $D_n$, when $n$ is even, we enumerate all Cartan subgroups (on the level of linear groups) as $\{   H _j ^d , 0\leqslant  j \leqslant n \}$, where the real rank of $H_j ^d$ is $j$, and we use the superscript $d$ to distinguish Cartan subgroups of the same real rank but not conjugate to each other. For example, when $n=4$, there are three Cartan subgroups of real rank 2, and they are labeled by $H_2 ^1$, $H_2 ^2,H_2 ^3$, all of which are isomorphic to $\mathbb{R}^{\times}\times S^1 \times \mathbb{C}^{\times}$.

 Let $\gamma _n$ be the parameter of the principal series. We start with a set of orthogonal nonintegral real roots $R(\gamma _n) = \{  \alpha _k, \beta _k , 1 \leqslant k \leqslant n/2\}$ of $\gamma _n$, where $\alpha _k = e_{2k-1}-e_{2k}$, $\beta _k = e_{2k-1} + e_{2k}$, and obtain $\gamma _j ^d$ by taking Cayley transforms through the roots in $R(\gamma _n)$. We attach to each $\gamma _j ^d$ a set of real roots $R(\gamma _j ^d)  = \{  \beta \in R(\gamma _n) \thinspace  |\thinspace   \beta \text{ is real for $\gamma _j ^d$ } \}$.  Now let $\gamma_0$ be the parameter of the discrete series with $R(\gamma _0)=\phi$, $\gamma _2 ^1$ be the parameter two steps up from $\gamma _0$, with $R(\gamma _2 ^1) =\{\alpha _{n/2}, \beta _{n/2}\}$, $\gamma ^2 _{n/2}$ be the parameter with $R(\gamma _{n/2} ^2) = \{\beta _1, \cdots , \beta _{n/2}\}$, and $\gamma ^3 _{n/2}$ be the parameter with $R(\gamma _{n/2} ^3) = \{\beta _1, \cdots , \beta _{n/2-1}, \alpha _{n/2} \}$. Observe that when $n=4$, $\gamma _2 ^1$ is the  representative from the middle Cartan subgroup $H^1 _2$; when $n>4$, choose $\gamma ^1 _{n/2}$ to be the representative from $H_{n/2} ^1$ with $R(\gamma _{n/2} ^1) = \{   \beta _2, \cdots, \beta _{n/2}, \alpha _{n/2}\}$. Note that it is possible that there exists $\gamma _{n/2} ^2$, $d>3$. In this case, choose $\gamma _{n/2} ^d$ such that $\{ \alpha _{n/2-1}, \alpha_{n/2}, \beta _{n/2-1},\beta_{n/2}\} \subseteq R(\gamma_{n/2} ^d )$. 

Now we claim that if $\gamma _j ^d$ is not one of these four, then it satisfies either condition (R) or (C), and hence can be ruled out by Lemma \ref{Druleout}.

When $j \geqslant n/2 +2$, $\gamma _j ^d$ can be chosen so that $\{\alpha_{n/2 -1}, \alpha _{n/2}, \beta _{n/2-1}, \beta _{n/2}\} \subseteq R(\gamma _j ^d)$ and hence $e_{n-2}\pm e_n$ are real integral roots of $\gamma _j ^d$.

Now suppose $n>4$. We observe that  $\gamma _{n/2} ^1$ satisfies condition (C) since there is a a quadruple $\{\alpha_{n/2 -1}, \alpha _{n/2}, \beta _{n/2-1}, \beta _{n/2}\}$, where $\alpha _{n/2}, \beta_{n/2}, \beta _{n/2-1}$ are real, and $\alpha _{n/2-1}$ is imaginary for $\gamma _{n/2} ^1$. For $d>3$,  since $\{ \alpha _{n/2-1}, \alpha_{n/2}, \beta _{n/2-1},\beta_{n/2}\} \subseteq R(\gamma_{n/2} ^d )$, $e_{n-2}\pm e_n$ are real imaginary roots of $\gamma _{n/2} ^d$, and hence $\gamma_{n/2} ^d$  satisfies condition (R). 

Any $\gamma ^d _{n/2 +1}$ is obtained from some $\gamma ^{d'} _{n/2}$ by an inverse Cayley transform, that is, we can choose $\gamma _{n/2 +1} ^d$ such that $R(\gamma _{n/2 +1} ^d )$ is obtained from $R(\gamma ^{d'} _{n/2})$ by adding a root. But adding a root to $R(\gamma _{n/2} ^{d'})$ would result in either a real integral roots or a quadruple as  described in condition (C)  for $\gamma _{n/2 +1} ^d$. 

Finally we observe $\gamma _{j} ^d$, $j> n/2$. Every $\gamma _{n/2} ^d$ can be obtained from some $\gamma _{n/2} ^{d'}$ by a sequence of Cayley transforms through roots in $R(\gamma _{n/2} ^{d'})$, that is,  $\gamma _j ^{d}$ is chosen such that $R_{\gamma _j} ^d$ is obtained by removing roots from $R(\gamma _{n/2} ^{d'})$. It turns out that when $j>n/2$, there would be a quadruple as described in condition (C) for all $\gamma _j ^d$, except $\gamma _0$ and $\gamma _2 ^1$. We conclude that $\gamma _0$, $\gamma _2 ^1$, $\gamma _ {n/2} ^2$, $\gamma _{n/2} ^3$ are the $\gamma _j$'s satisfying condition $(\ast)$ since the number of $\prod _{R_D} (  \widetilde{G})$ with a fixed central character is exactly 4. Hence the theorem follows for type $D_n$, $n \geqslant 4$, when $n$ is even. 

When $n$ is odd, it can be shown by a similar argument that the two $\gamma_j$'s satisfying condition $(\ast)$ come from the compact Cartan subgroup and a Cartan subgroup with real rank 2. 
\end{proof}
\end{theorem}
We would like to do the same thing for type $E$, parallel to the case of type $D_n$. Like in type $D_n$, we can also move to a block $\mathcal{D}$ where all integral simple roots are simple but one, say $\alpha$. Then there come some difficulties. First, to decompose $s_{\alpha}$ into a product of simple reflections is never easy, and after having done  that, we have to keep track of a sequence of inner automorphisms when trying to calculate the coherent continuation action $s_{\alpha}\cdot \gamma$, where $\gamma$ is a standard module parameter. Even though this complication has not been solved yet, we strongly believe that the counting $ |\prod _{\rho /2} ^s (\widetilde{G})|=|  \prod _{R_D} (\widetilde{G}) |$ in Theorem \ref{exhaustAD} holds for type $E_6, E_7$ and $E_8$.

We conjecture that when we count the number of $\gamma_j$ satisfying condition $(\ast)$, those satisfying condition (R) or (C) should be ruled out  and hence Theorem \ref{exhaustAD} is true for type $E$, that is,  Conjecture \ref{conjE} would be true if Conjecture \ref{Eruleout} is true.

\begin{conjecture} \label{Eruleout}
For type $E_6$, $E_7$ and $E_8$, $\gamma _j$ does not satisfy condition $(\ast)$ if it satisfy condition (R) or (C).
\end{conjecture}

\begin{conjecture}\label{conjE}
For type $E_6$, $E_7$ and $E_8$, we have  $\prod _{\lambda} ^s (\widetilde{G})= \prod _{R_D}(\widetilde{G})$.
\end{conjecture}

\section{Relation to the pairs $(\chi, \mathcal{O}_{\mathbb{R}})$ }
In this section, the setting is the same as in Section \ref{splitcase}, that is, $G$ is simply laced and split.  We define
\begin{center}
$P_{\Or} (\tu{G})=\{  (\chi _i, \mathcal{O} _j ) \thinspace| \thinspace   \chi _i \in \prod _g (Z(\tu{G})), \mathcal{O}_j \text{ is a real form of } \mathcal{O} \}$.
\end{center}

 Due to Proposition \ref{singlereal}, there is a mapping from $\tu{\pi}\in \prod _{\rho/2} ^s (\tu{G} )$ to  $P_{\Or} (\tu{G})$ given by 
 \begin{equation}\label{corresp}
\xi: \tu{\pi} \longmapsto (\chi _{\tu{\pi}}, \text{AV}(\tu{\pi}) ),
 \end{equation}
where $\chi _{\tu{\pi}}$ is the central character of $\tu{\pi}$. By Proposition \ref{Shi}(3) and Remark \ref{number},  $|P_{\Or} (\tu{G})| = | \prod_{R_D} (\tu{G})|$, and this number is equal to $|\prod_{\rho/2} ^s (\tu{G}) | $ when $G$ has type $A_{n-1}$ or $D_n$. We will take a look at the correspondence between $\prod_{\rho/2} ^s (\tu{G})$ and $P_{\Or} (\tu{G})$  in the cases of type $A$ and $D$.

\subsection{Type $A_{n-1}$}
When $n$ is odd,  $\Smal =\SH _{\rho/2}$, which contains a single  $Sh_i$, and hence $|P_{\Or} (\tu{G})| =1$. The mapping $\xi$ in (\ref{corresp}) is nothing but a one-to-one correspondence.

When $n$ is even, say, $n=2m$, $\Smal =\{ Sh_1, Sh_2, \pi_1, \pi _2\}$, where $\pi_1$ and $\pi _2$ are constructed in Section \ref{piRD} and representations with the same subscript have the same central character. The $K$-types of these representations are calculated in \cite{Luc}, and we list them as follows. By convention, the $K$-types of $Sh_1$ and $Sh_2$ are 
 \begin{center}$(\frac{1}{2}+2a_1, \cdots, \frac{1}{2}+2a_{m} )$, $a_1\geq a_2\geq \cdots \geq a_m\geq 0$,  $a_i\in\mathbb{Z}$ and \\ 
$(\frac{1}{2}+2a_1, \cdots, \frac{1}{2}+2a_{m-1}, -(\frac{1}{2}+2a_m) )$, $a_1\geq a_2\geq \cdots \geq a_m\geq 0$,  $a_i\in\mathbb{Z}$,
\end{center}
respectively. Because $\pi_i$ and $Sh_i$ have the same central character, when $m$ is even, the $K$-types of $\pi _1$ and $\pi_2$ are 
\begin{center}
$(\frac{3}{2}+2a_1, \cdots, \frac{3}{2}+2a_{m} )$, $a_1\geq a_2\geq \cdots \geq a_m\geq 0$,  $a_i\in\mathbb{Z}$ and \\ 
$(\frac{3}{2}+2a_1, \cdots, \frac{3}{2}+2a_{m-1}, -(\frac{3}{2}+2a_m) )$, $a_1\geq a_2\geq \cdots \geq a_m\geq 0$,  $a_i\in\mathbb{Z}$,
\end{center}
respectively; when $m$ is odd, the $K$-types of $\pi _1$ and $\pi_2$ are
\begin{center}
$(\frac{3}{2}+2a_1, \cdots, \frac{3}{2}+2a_{m-1}, -(\frac{3}{2}+2a_m) )$, $a_1\geq a_2\geq \cdots \geq a_m\geq 0$,  $a_i\in\mathbb{Z}$ and\\
$(\frac{3}{2}+2a_1, \cdots, \frac{3}{2}+2a_{m} )$, $a_1\geq a_2\geq \cdots \geq a_m\geq 0$,  $a_i\in\mathbb{Z}$,
\end{center}
respectively. Observe that $Sh_i$ and $\pi _i$, $i=1,2$, have the same asymptotic $K$-types when $m$ is odd; $Sh_1$ and $\pi_2$ have the same   asymptotic $K$-types and $Sh_2$ and $\pi_1$ have the same   asymptotic $K$-types when $m$ is even. By the convention made in Remark \ref{Shiorb}, we have  Table 5 and 6 to illustrate the correspondence $\xi$.

\begin{table}[H] \label{Aorbit}
\caption{$A_{2m-1}$,  $m$ even}
\begin{center}
\begin{tabular}
{c|c|c}&$\mathcal{O}_1 $& $\mathcal{O}_2 $ \\\hline $\chi _1$ &$Sh _1$, $\pi_1$ & N/A  \\  && \\ \hline $\chi_2$ & N/A &$Sh_2$, $\pi_2$ \\ && \end{tabular}
\end{center}
\end{table}
\begin{table}[H] 
\caption{$A_{2m-1}$,  $m$ odd}
\begin{center}
\begin{tabular}
{c|c|c}&$\mathcal{O}_1 $& $\mathcal{O}_2 $ \\\hline $\chi _1$ &$Sh _1$ & $\pi_1$  \\  &\footnotesize{$(\frac{1}{2},\cdots,\frac{1}{2})$}& \footnotesize{$(\frac{3}{2},\cdots,\frac{3}{2}, -\frac{3}{2})$}\\ \hline $\chi_2$ & $\pi_2$ &$Sh_2$ \\ & $\footnotesize{ (\frac{3}{2},\cdots,\frac{3}{2} ) }$ & \footnotesize{ $(\frac{1}{2},\cdots,\frac{1}{2}, -\frac{1}{2})$}
\end{tabular}
\end{center}
\end{table}

We conclude that $\xi$ is a bijection if and only if $n\equiv 1, 2, 3$ (mod 4).
 
\subsection{Type $D_n$}
In this case, the double cover $\tu{G} =\tu{\text{ Spin}}(n,n)$. In \cite{LS1}, there are some small representations  defined for  a bigger group $\widetilde{G'}=\widetilde{\text{Spin}}(n+1,n)$ with maximal compact subgroup $\tu{K'} =\text{Spin}(n+1)\times \text{Spin} (n)$. When $n=2m$, there are four of them, denoted $\Gamma_1, \Gamma_2, \Gamma_3, \Gamma_4$; when $n=2m+1$, there are two of them, denoted $\Gamma_1, \Gamma _2$. Their $K$-types are calculated in \cite{LS1}. For example, looking at Table 9, $Sh_3$ is obtained by restriction of $\Gamma_1$ to the subgroup $Spin(n,n)\times O(2)$ with $O(2)$ acting trivially; on the other hand, $\pi _4$  is obtained by restriction of $\Gamma_1$ to the subgroup $Spin(n,n)\times O(2)$ with $O(2)$ acting by the sign character.

Let  $\tu{K}=\text{Spin}(2n)\times \text{Spin} (2n)$, the maximal compact subgroup of $\widetilde{G}$. Suppose $n=2m$. Then  Out$(\tu{G})$ is generated by two elements $\sigma$ and $\gamma$ and the action of them on $K$-types are as follows. \begin{center} $\sigma(\lambda_1,\cdots,\lambda_m; \lambda_{m+1} ,\cdots, \lambda_{n})=  (\lambda_1,\cdots,-\lambda_m; \lambda_{m+1} ,\cdots, -\lambda_{n})$, \\$\gamma(\lambda_1,\cdots,\lambda_m; \lambda_{m+1} ,\cdots, \lambda_{n})=(\lambda_{m+1} ,\cdots, \lambda_{n}  ; \lambda_1,\cdots,\lambda_m)$ \end{center}
The $K$-types  parametrized by $(\lambda ;\lambda ')$ and $(\lambda' ; \lambda)$ represent different representations when restricted to $\tu{K}$ since $\gamma$ is in Out($\tu{G}$), and hence  restricting these $\Gamma_i$'s to $\widetilde{G}$, we  get 16 representations (see Table 9 and 10 for the labeling for these representations and their $K$-types), and they are the small representations in $\Smal$. Notice that in Table 9 and 10,  $\lambda=(\lambda _1, \cdots , \lambda _n) \in \mathbb{Z}^n$ and  $\gamma =(\gamma_1,\cdots,\gamma_n)\in \mathbb{Z}^n$ are decreasing sequences of nonnegative integers. By $\gamma \prec\lambda$ we mean that $\lambda_1 \geqslant \gamma_1 \geqslant \cdots \geqslant \lambda _n\geqslant \gamma_n \geqslant -\lambda _n$. Also, $\bold{0}=(0,\cdots, 0), \bold{1}=(1,\cdots,1), \bold{\frac{1}{2}} =(\frac{1}{2},\cdots, \frac{1}{2})$. 
 
Suppose $n=2m+1$. Then Out($\tu{G}$) is generated by $\sigma$ (defined as above). Restricting $\Gamma _i$'s to $\tu{G}$ gives four representations and they are the ones in $\Smal$ (see Table 11 for their $K$-types). Notice that in Table 11, $\lambda=(\lambda _1, \cdots , \lambda _n) \in \mathbb{Z}^n$, $\lambda '=(\lambda _1, \cdots , \lambda _n, 0) \in \mathbb{Z}^{n+1}$ and  $\gamma =(\gamma_1,\cdots,\gamma_n)\in \mathbb{Z}^n$ are decreasing sequences of nonnegative integers. By $\gamma \prec\lambda '$ we mean that $\lambda_1 \geqslant \gamma_1 \geqslant \cdots \geqslant \lambda _n\geqslant \gamma_n \geqslant 0$. 

Observing the asymptotic $K$-types of the representations in Table 9, 10 and 11, it is illustrated in Table 7 and 8 that  $\xi$ is a bijection.

The following theorem follows immediately from the above discussion.
\begin{theorem}\label{sp}
Assume the setting in this section. If $G$ has type $A_{n-1}$, n is not divisible by 4, or $D_{2m}$, then $\xi$ gives a one-to-one correspondence between $\Smal$ and $P_{\Or} (\tu{G})$. 
\end{theorem}

Notice that for type $E_6$ and $E_8$, $\xi$ is a one-to-one correspondence between $\prod _{R_D}(\tu{G})$ and $P_{\Or} (\tu{G})$ since both of them contain a single  $Shi$, and hence $\xi$ is surjective map from $\Smal$ to $P_{\Or}(\tu{G} ).$ Because of Conjecture \ref{conjE}, we conjecture that Theorem \ref{sp} is true for type $E$. 
\begin{conjecture}
Assume the setting in this section. If $G$ has type $E_6, E_7$ or $E_8$, then $\xi$ is a one-to-one correspondence between $\Smal$ and $P_{\Or} (\tu{G})$. 
\end{conjecture}
\begin{table}[H]
\caption{$D_{2m}$}
\begin{center}
\begin{tabular}{|c|c|c|c|c|}\hline  &$\mathcal{O}_1 $& $\mathcal{O}_2 $ & $\mathcal{O}_3 $& $\mathcal{O}_4$ \\\hline $\chi _1$ &$Sh _1$ & $\pi_1$ & $\delta _1$ & $\tau _1$ \\  &\footnotesize{$(\frac{1}{2},\cdots,\frac{1}{2};0,\cdots,0)$}& \footnotesize{$(\frac{3}{2},\frac{1}{2},\cdots,-\frac{1}{2};0,\cdots,0)$}&\footnotesize{$(\frac{1}{2},\cdots,\frac{1}{2};1,\cdots,1)$}& \footnotesize{$(\frac{1}{2},\cdots,\frac{1}{2};1,\cdots,-1)$}\\ \hline $\chi_2$ & $\pi_2$ &$Sh_2$ & $\tau _2$& $\delta _2$\\ & \footnotesize{$(\frac{3}{2},\frac{1}{2},\cdots,\frac{1}{2};0,\cdots,0)$}& \footnotesize{$(\frac{1}{2},\cdots,-\frac{1}{2};0,\cdots,0)$}&\footnotesize{$(\frac{1}{2},\cdots,-\frac{1}{2};1,\cdots,1)$} &\footnotesize{ $(\frac{1}{2},\cdots,-\frac{1}{2};1,\cdots,-1)$} \\\hline $\chi _3$ & $\delta _3$ & $\tau_3$ & $Sh_3$ & $\pi _3$ \\ &\footnotesize{$(1,\cdots,1; \frac{1}{2},\cdots,\frac{1}{2})$}&\footnotesize{$(1,\cdots,-1; \frac{1}{2},\cdots,\frac{1}{2})$}&\footnotesize{ $(0,\cdots,0;\frac{1}{2},\cdots,\frac{1}{2})$}& \footnotesize{$(0,\cdots,0;\frac{3}{2},\frac{1}{2},\cdots,-\frac{1}{2})$ } \\\hline $\chi _4$ & $\tau_4$& $\delta_4$& $\pi _4$& $Sh_4$ \\ &\footnotesize{$(1,\cdots,1; \frac{1}{2},\cdots,-\frac{1}{2})$}&\footnotesize{$(1,\cdots,-1; \frac{1}{2},\cdots,-\frac{1}{2})$} & \footnotesize{$(0,\cdots,0;\frac{3}{2},\frac{1}{2},\cdots,\frac{1}{2})$}& \footnotesize{$(0,\cdots,0;\frac{1}{2},\cdots,-\frac{1}{2})$}\\\hline \end{tabular}
\end{center}
\end{table}

\begin{table}[H]
\caption{$D_{2m+1}$}
\begin{center}
\begin{tabular}{|c|c|c|c|c|}\hline  &$\mathcal{O}_1 $& $\mathcal{O}_2 $ \\\hline $\chi _1$ &$Sh _1$,  $\pi_1$&  N/A \\  &&\\ \hline $\chi_2$ &N/A&  $Sh_2$, $\pi_2$ \\ &  &\\\hline \end{tabular}
\end{center}
\end{table}

\begin{table}[H] \label{DKtype1}
\caption{All $K$-types of representations in $\Smal$ when $G=Spin(n,n)$, $n=2m$ (1)}  
\scriptsize
\begin{center}
\begin{tabular}{|c|c|c|c|c|c|}\hline $\widetilde{G'}$-rep &  $K'$-type& L.$K'$.T & restriction to $K$& L.$K$.T & $\widetilde{G}$-rep  \\\hline \multirow{2}{*}{ $\Gamma _1  $} & \multirow{2}{*}{ $ V^+ =\bigoplus  \limits _{\lambda}(\lambda; \lambda +\displaystyle\bold{ \frac{1}{2}})$}&  \multirow{2}{*}{ $(\bold{0};\bold{\displaystyle \frac{1}{2}})$} & $\bigoplus \limits _{\lambda}\bigoplus  \limits _{\gamma \prec \lambda} (\gamma; \lambda +\displaystyle\bold{ \frac{1}{2}}),  \Sigma (\lambda _i+\gamma _i)\in 2\mathbb{Z}$  &  $(0,\cdots,0;\frac{1}{2},\cdots,\frac{1}{2})$& $Sh _3$ \\ \cline{4-4}\cline{5-5} \cline{6-6}  &   & & $\bigoplus \limits _{\lambda}\bigoplus  \limits _{\gamma \prec \lambda }(\gamma; \lambda +\displaystyle\bold{ \frac{1}{2}}),  \Sigma (\lambda _i+\gamma _i)\in 2\mathbb{Z}+1$ &$(0,\cdots,0;\frac{3}{2},\frac{1}{2},\cdots,\frac{1}{2})$ & $\pi _4$   \\\hline  \multirow{2}{*}{  $ \Gamma _1$} & \multirow{2}{*}{ $ V^+ =\bigoplus  \limits _{\lambda} (\lambda +\displaystyle\bold{ \frac{1}{2}} ; \lambda)$}&  \multirow{2}{*}{ $(\bold{\displaystyle \frac{1}{2}} ; \bold{0})$} & $\bigoplus \limits _{\lambda}\bigoplus  \limits _{\gamma \prec \lambda }(\lambda +\displaystyle\bold{ \frac{1}{2}};\gamma),  \Sigma (\lambda _i+\gamma _i)\in 2\mathbb{Z}$  &  $(\frac{1}{2},\cdots,\frac{1}{2};0,\cdots,0)$& $Sh _1$ \\ \cline{4-4}\cline{5-5} \cline{6-6}  &   & &  $\bigoplus \limits _{\lambda}\bigoplus  \limits _{\gamma \prec \lambda }(\lambda +\displaystyle\bold{ \frac{1}{2}}; \gamma),  \Sigma (\lambda _i+\gamma _i)\in 2\mathbb{Z}+1$ &$(\frac{3}{2},\frac{1}{2},\cdots,\frac{1}{2};0,\cdots,0)$ & $\pi _2$  \\\hline   \multirow{2}{*}{$\Gamma _2 $  } & \multirow{2}{*}{ $ V^- =\bigoplus  \limits _{\lambda}(\lambda; \sigma(\lambda +\displaystyle\bold{ \frac{1}{2}}))$}&  \multirow{2}{*}{ $(\bold{0};\sigma(\bold{\displaystyle \frac{1}{2}}))$} & $\bigoplus \limits _{\lambda}\bigoplus  \limits _{\gamma\prec\lambda }(\gamma; \sigma(\lambda +\displaystyle\bold{ \frac{1}{2}})),  \Sigma (\lambda _i+\gamma _i)\in 2\mathbb{Z}$  &  $(0,\cdots,0;\frac{1}{2},\cdots,-\frac{1}{2})$& $Sh_4$ \\ \cline{4-4}\cline{5-5} \cline{6-6}  &   & & $\bigoplus \limits _{\lambda}\bigoplus  \limits _{\gamma \prec \lambda }(\gamma; \sigma(\lambda +\displaystyle\bold{ \frac{1}{2}})),  \Sigma (\lambda _i+\gamma _i)\in 2\mathbb{Z}+1$ &$(0,\cdots,0;\frac{3}{2},\frac{1}{2},\cdots,-\frac{1}{2})$ &   $\pi _3$\\\hline  \multirow{2}{*}{  $\Gamma _2$} & \multirow{2}{*}{ $ V^- =\bigoplus  \limits _{\lambda} (\sigma(\lambda +\displaystyle\bold{ \frac{1}{2}}) ; \lambda)$}&  \multirow{2}{*}{ $(\sigma(\bold{\displaystyle \frac{1}{2}}) ; \bold{0})$} & $\bigoplus \limits _{\lambda}\bigoplus  \limits _{\gamma \prec\lambda }(\sigma(\lambda +\displaystyle\bold{ \frac{1}{2}}); \gamma),  \Sigma (\lambda _i+\gamma _i)\in 2\mathbb{Z}$  &  $(\frac{1}{2},\cdots,-\frac{1}{2};0,\cdots,0)$&$Sh _2$ \\ \cline{4-4}\cline{5-5} \cline{6-6}  &   & &  $\bigoplus \limits _{\lambda}\bigoplus  \limits _{\gamma \prec \lambda }(\sigma(\lambda +\displaystyle\bold{ \frac{1}{2}}); \gamma),  \Sigma (\lambda _i+\gamma _i)\in 2\mathbb{Z}+1$ &$(\frac{3}{2},\frac{1}{2},\cdots,-\frac{1}{2};0,\cdots,0)$ & $\pi _1$  \\\hline 
 
\end{tabular}
\end{center}
\end{table}

\begin{table} \label{DKtype2}
\caption{All $K$-types of representations in $\Smal$ when $G=Spin(n,n)$, $n=2m$ (2)}  
\scriptsize
\begin{center}
\begin{sideways}
\begin{tabular}{|c|c|c|c|c|c|} 
\hline $\widetilde{G'}$-rep &  $K'$-type& L.$K'$.T & restriction to $K$& L.$K$.T & $\widetilde{G}$-rep  \\
\hline
 \multirow{4}{*}{$ \Gamma_3 $} & \multirow{4}{*}{$ V_0 ^+ =\bigoplus  \limits _{\lambda} (\lambda +\displaystyle\bold{ \frac{1}{2}} ; \lambda+\bold{1})$}&  \multirow{4}{*}{ $(\bold{\displaystyle \frac{1}{2}} ; \bold{1})$} &   $\bigoplus  \limits _{\lambda} \bigoplus  \limits _{\gamma \prec \lambda }(\gamma +\displaystyle\bold{ \frac{1}{2}}; \lambda+\bold{1}),  \Sigma (\lambda _i+\gamma _i)\in 2\mathbb{Z}$ &  \multirow{2}{*}{$(\frac{1}{2},\cdots,\frac{1}{2};1,\cdots,1)$}& \multirow{2}{*}{$\delta _1 (m \text{ even})/ \delta _2 (m \text{ odd})$} \\ &   & &  $\bigoplus  \limits _{\lambda} \bigoplus  \limits _{\gamma \prec \lambda }(\sigma (\gamma +\displaystyle\bold{ \frac{1}{2}}); \lambda+\bold{1}),  \Sigma (\lambda _i+\gamma _i)\in 2\mathbb{Z}+1$& & \\ \cline{4-4}\cline{5-5}\cline{6-6} &&& 
 $\bigoplus  \limits _{\lambda} \bigoplus  \limits _{\gamma \prec \lambda }(\gamma +\displaystyle\bold{ \frac{1}{2}}; \lambda+\bold{1}),  \Sigma (\lambda _i+\gamma _i)\in 2\mathbb{Z}+1$ &\multirow{2}{*}{ $(\frac{1}{2},\cdots,\frac{1}{2}, -\frac{1}{2};1,\cdots,1)$} & \multirow{2}{*}{ $\tau _2 (m \text{ even})/ \tau _1 (m \text{ odd})$}\\ &&& $\bigoplus  \limits _{\lambda} \bigoplus  \limits _{\gamma \prec \lambda }(\sigma (\gamma +\displaystyle\bold{ \frac{1}{2}}); \lambda+\bold{1},  \Sigma (\lambda _i+\gamma _i)\in 2\mathbb{Z}$ &&  \\\hline 
 \multirow{4}{*}{$ \Gamma_3 $} & \multirow{4}{*}{ $ V_0 ^+ =\bigoplus  \limits _{\lambda} (\lambda+\bold{1}; \lambda +\displaystyle\bold{ \frac{1}{2}} )$}&  \multirow{4}{*}{ $(\bold{1}; \bold{\displaystyle \frac{1}{2}} )$} & $\bigoplus \limits_{\lambda} \bigoplus  \limits _{\gamma \prec \lambda }(\lambda+\bold{1} ; \gamma+\displaystyle\bold{ \frac{1}{2}}),  \Sigma (\lambda _i+\gamma _i)\in 2\mathbb{Z}$  &  \multirow{2}{*}{ $(1,\cdots,1; \frac{1}{2},\cdots,\frac{1}{2})$}& \multirow{2}{*}{$\delta _3 (m \text{ even})/ \delta _4 (m \text{ odd})$} \\ &   & &  $\bigoplus \limits_{\lambda} \bigoplus  \limits _{\gamma \prec \lambda }(\lambda+\bold{1}; \sigma(\gamma+\displaystyle\bold{ \frac{1}{2}}     )  ),  \Sigma (\lambda _i+\gamma _i)\in 2\mathbb{Z}+1$ & & \\ \cline{4-4}\cline{5-5}\cline{6-6} &&& 
  $\bigoplus \limits_{\lambda}  \bigoplus  \limits _{\gamma \prec \lambda } (\lambda + \bold{1}; \gamma +\displaystyle\bold{ \frac{1}{2}}),  \Sigma (\lambda _i+\gamma _i)\in 2\mathbb{Z}+1$ &\multirow{2}{*}{
  $(1,\cdots,1; \frac{1}{2},\cdots,-\frac{1}{2})$} & \multirow{2}{*}{ $\tau _4 (m \text{ even})/ \tau _3 (m \text{ odd})$}\\ &&& $\bigoplus \limits_{\lambda} \bigoplus  \limits _{\gamma \prec \lambda }(\lambda+\bold{1} ; \sigma(\gamma+\displaystyle\bold{ \frac{1}{2}})),  \Sigma (\lambda _i+\gamma _i)\in 2\mathbb{Z}$  &&  \\\hline 
 \multirow{4}{*}{$ \Gamma_4 $} & \multirow{4}{*}{ $ V_0 ^- =\bigoplus  \limits _{\lambda} (\lambda +\displaystyle\bold{ \frac{1}{2}} ; \sigma(\lambda+\bold{1}))$}&  \multirow{4}{*}{ $(\bold{\displaystyle \frac{1}{2}} ; \sigma(\bold{1}))$} &   $\bigoplus  \limits _{\lambda} \bigoplus  \limits _{\gamma \prec \lambda }(\gamma +\displaystyle\bold{ \frac{1}{2}}; \sigma(\lambda+\bold{1})),  \Sigma (\lambda _i+\gamma _i)\in 2\mathbb{Z}$ &  \multirow{2}{*}{$(\frac{1}{2},\cdots,\frac{1}{2};1,\cdots,-1)$}& \multirow{2}{*}{$\tau _1 (m \text{ even})/ \tau _2 (m \text{ odd})$} \\ &   & &  $\bigoplus  \limits _{\lambda} \bigoplus  \limits _{\gamma \prec \lambda }(\sigma (\gamma +\displaystyle\bold{ \frac{1}{2}}); \sigma(\lambda+\bold{1})),  \Sigma (\lambda _i+\gamma _i)\in 2\mathbb{Z}+1$& & \\ \cline{4-4}\cline{5-5}\cline{6-6} &&& 
 $\bigoplus  \limits _{\lambda} \bigoplus  \limits _{\gamma \prec \lambda }(\gamma +\displaystyle\bold{ \frac{1}{2}}; \sigma(\lambda+\bold{1})),  \Sigma (\lambda _i+\gamma _i)\in 2\mathbb{Z}+1$ &\multirow{2}{*}{
 $(\frac{1}{2},\cdots,-\frac{1}{2};1,\cdots,-1)$} & \multirow{2}{*}{ $\delta _2 (m \text{ even})/ \delta _1 (m \text{ odd})$}\\ &&& $\bigoplus  \limits _{\lambda} \bigoplus  \limits _{\gamma \prec \lambda }(\sigma( \gamma +\displaystyle\bold{ \frac{1}{2}}; \sigma(\lambda+\bold{1})),  \Sigma (\lambda _i+\gamma _i)\in 2\mathbb{Z}$ &&  \\\hline 

 \multirow{4}{*}{$ \Gamma_4 $} & \multirow{4}{*}{ $ V_0 ^- =\bigoplus  \limits _{\lambda} (\sigma(\lambda+\bold{1}); \lambda +\displaystyle\bold{ \frac{1}{2}} )$}&  \multirow{4}{*}{ $(\sigma(\bold{1}); \bold{\displaystyle \frac{1}{2}} )$} & $\bigoplus \limits_{\lambda} \bigoplus  \limits _{\gamma \prec \lambda }(\sigma(\lambda+\bold{1}) ; \gamma+\displaystyle\bold{ \frac{1}{2}}),  \Sigma (\lambda _i+\gamma _i)\in 2\mathbb{Z}$  &  \multirow{2}{*}{ $(1,\cdots,-1; \frac{1}{2},\cdots,\frac{1}{2})$}& \multirow{2}{*}{$\tau _3 (m \text{ even})/ \tau _4 (m \text{ odd})$} \\ &   & &  $\bigoplus \limits_{\lambda} \bigoplus  \limits _{\gamma \prec \lambda }(\sigma(\lambda+\bold{1}) ; \sigma(\gamma+\displaystyle\bold{ \frac{1}{2}}     )  ),  \Sigma (\lambda _i+\gamma _i)\in 2\mathbb{Z}+1$ & & \\ \cline{4-4}\cline{5-5}\cline{6-6} &&& 
  $\bigoplus \limits_{\lambda}  \bigoplus  \limits _{\gamma \prec \lambda } (\sigma(\lambda + \bold{1}); \gamma +\displaystyle\bold{ \frac{1}{2}}),  \Sigma (\lambda _i+\gamma _i)\in 2\mathbb{Z}+1$ &\multirow{2}{*}{
  $(1,\cdots,-1; \frac{1}{2},\cdots,-\frac{1}{2})$} & \multirow{2}{*}{ $\delta _4 (m \text{ even})/ \delta _3 (m \text{ odd})$}\\ &&& $\bigoplus \limits_{\lambda} \bigoplus  \limits _{\gamma \prec \lambda }(\sigma(\lambda+\bold{1}) ; \sigma(\gamma+\displaystyle\bold{ \frac{1}{2}})),  \Sigma (\lambda _i+\gamma _i)\in 2\mathbb{Z}$  &&  \\\hline 

\end{tabular}
\end{sideways}
\end{center}
\end{table}
\begin{table} \label{DKtype3}
\caption{All $K$-types of representations in $\Smal$, when $G=Spin(n,n)$, $n=2m+1$ }
\scriptsize
\begin{center}
%\begin{sideways}
\begin{tabular}{|c|c|c|c|c|c|}\hline $\widetilde{G'}$-rep &  $K'$-type& L.$K'$.T & restriction to $K$& L.$K$.T & $\widetilde{G}$-rep  \\\hline \multirow{2}{*}{ $\Gamma   $} & \multirow{2}{*}{ $ V =\bigoplus  \limits _{\lambda}(\lambda '; \lambda +\displaystyle\bold{ \frac{1}{2}})$}&  \multirow{2}{*}{ $(\bold{0};\bold{\displaystyle \frac{1}{2}})$} & $\bigoplus \limits _{\lambda}\bigoplus  \limits _{\gamma \prec \lambda '} (\gamma; \lambda +\displaystyle\bold{ \frac{1}{2}}),  \Sigma (\lambda _i+\gamma _i)\in 2\mathbb{Z}$  &  $(0,\cdots,0;\frac{1}{2},\cdots,\frac{1}{2})$& $Sh _2$ \\ \cline{4-4}\cline{5-5} \cline{6-6}  &   & & $\bigoplus \limits _{\lambda}\bigoplus  \limits _{\gamma \prec \lambda ' }(\gamma; \lambda +\displaystyle\bold{ \frac{1}{2}}),  \Sigma (\lambda _i+\gamma _i)\in 2\mathbb{Z}+1$ &$(0,\cdots,0;\frac{3}{2},\frac{1}{2},\cdots,\frac{1}{2})$ & $\pi _2$   \\\hline  \multirow{2}{*}{  $ \Gamma $} & \multirow{2}{*}{ $ V =\bigoplus  \limits _{\lambda} (\lambda +\displaystyle\bold{ \frac{1}{2}} ; \lambda')$}&  \multirow{2}{*}{ $(\bold{\displaystyle \frac{1}{2}} ; \bold{0})$} & $\bigoplus \limits _{\lambda}\bigoplus  \limits _{\gamma \prec \lambda' }(\lambda +\displaystyle\bold{ \frac{1}{2}}; \gamma),  \Sigma (\lambda _i+\gamma _i)\in 2\mathbb{Z}$  &  $(\frac{1}{2},\cdots,\frac{1}{2};0,\cdots,0)$& $Sh _1$ \\ \cline{4-4}\cline{5-5} \cline{6-6}  &   & &  $\bigoplus \limits _{\lambda}\bigoplus  \limits _{\gamma \prec \lambda' }(\lambda +\displaystyle\bold{ \frac{1}{2}}; \gamma),  \Sigma (\lambda _i+\gamma _i)\in 2\mathbb{Z}+1$ &$(\frac{3}{2},\frac{1}{2},\cdots,\frac{1}{2};0,\cdots,0)$ & $\pi _1$  \\\hline 
\end{tabular}
%\end{sideways}
\end{center}
\end{table}
\newpage

\end{document}